\documentclass[]{interact}

\usepackage{epstopdf}
\usepackage[caption=false]{subfig}

\usepackage[numbers,sort&compress]{natbib}
\bibpunct[, ]{[}{]}{,}{n}{,}{,}
\renewcommand\bibfont{\fontsize{10}{12}\selectfont}

\theoremstyle{plain}
\newtheorem{theorem}{Theorem}[section]
\newtheorem{lemma}[theorem]{Lemma}
\newtheorem{corollary}[theorem]{Corollary}

\theoremstyle{definition}
\newtheorem{definition}[theorem]{Definition}

\theoremstyle{remark}
\newtheorem{remark}{Remark}

\usepackage{comment}
\usepackage{color}
\usepackage[colorlinks=true,breaklinks=true,bookmarks=true,urlcolor=blue,citecolor=blue,linkcolor=blue,bookmarksopen=false,draft=false]{hyperref}
\usepackage{algorithm}
\usepackage{algcompatible}
\usepackage{graphicx}
\usepackage{caption}
\usepackage{multirow,booktabs,dcolumn}
\usepackage[export]{adjustbox}
\usepackage[utf8]{inputenc}

\newcommand{\executeiffilenewer}[3]{%
	\ifnum\pdfstrcmp{\pdffilemoddate{#1}}%
	{\pdffilemoddate{#2}}>0%
	{\immediate\write18{#3}}\fi%
}

\graphicspath{{fig/}}
\begin{document}


\title{A gradient descent akin method for constrained optimization: algorithms and applications}

\author{
\name{Long Chen\textsuperscript{a}\thanks{L. Chen. Email: long.chen@scicomp.uni-kl.de}, Kai-Uwe Bletzinger\textsuperscript{b}\thanks{K.-U. Bletzinger. Email: kub@tum.de}, Nicolas R. Gauger\textsuperscript{a}\thanks{N. R. Gauger. Email: nicolas.gauger@scicomp.uni-kl.de} and Yinyu Ye\textsuperscript{c}\thanks{Y. Ye. Email: yyye@stanford.edu}}
\affil{\textsuperscript{a}Chair for Scientific Computing, University of Kaiserslautern-Landau (RPTU), Germany;\\ \textsuperscript{b}Chair of Structural Analysis, Technical University of Munich, Germany; \\
\textsuperscript{c}Department of Management Science and Engineering and ICME, Stanford University, USA}
}

\maketitle

\begin{abstract}
We present a first-order method for solving constrained optimization problems. The method is derived from our previous work \cite{Chen}, a modified search direction method inspired by singular value decomposition. In this work, we simplify its computational framework to a ``gradient descent akin'' method (GDAM), i.e., the search direction is computed using a linear combination of the negative and normalized objective and constraint gradient. We give fundamental theoretical guarantees on the global convergence of the method. This work focuses on the algorithms and applications of GDAM. We present computational algorithms that adapt common strategies for the gradient descent method. We demonstrate the potential of the method using two engineering applications, shape optimization and sensor network localization. When practically implemented, GDAM is robust and very competitive in solving the considered large and challenging optimization problems.
\end{abstract}


\begin{amscode}
65K05; 90C22; 90C30; 90C51; 90C90
\end{amscode}

%

\section{Introduction}
\label{sec:intro}
In this paper, we propose a first-order method for solving inequality constrained optimization problems. The problem of interest is
\begin{equation}
\tag{COP}
\begin{split}
&\textnormal{minimize}~~~  f(x), \\
&\textnormal{subject to}~~ g_i(x) \leq 0, ~ i = 1,...,m, \\
\end{split}
\label{eq:Optimization_Problem}
\end{equation}
where $f, g_1,..., g_m: \mathbb{R}^{n} \rightarrow \mathbb{R}$ are twice differentiable. We propose a gradient descent akin method (GDAM) for \eqref{eq:Optimization_Problem} with a single inequality constraint $g(x)$,
\begin{equation}
\tag{GDAM}
\begin{split}
    &x_{k+1} = x_k + \alpha_k \mathbf{s}_\zeta (x_k), \\
    &\mathbf{s}_\zeta (x_k) = -\frac{\nabla f(x_k)}{|\nabla f(x_k)|}
-\zeta\frac{\nabla g(x_k)}{|\nabla g(x_k)|}, ~~ \zeta \in [0,1),
\end{split}  
\label{eq:GDAM}
\end{equation}
where $\alpha_k$ denotes the step size, ``$|\cdot|$" denotes the Euclidean norm, and $\zeta$ is a parameter. $\nabla f(x)$ and $\nabla g(x)$ are the gradient column vectors of the objective function and constraint function, respectively. 

Provided the well-known logarithmic barrier function
\begin{equation}
\Phi(x) = - \sum_{i = 1}^{m} \log (-g_i(x)), ~ i = 1,...,m,
\label{eq:log_barrier}
\end{equation}
where $\log(\cdot)$ denotes the natural logarithm, 
we propose a generalization of \eqref{eq:GDAM} for multiple constrained problems with $\zeta\in [0,1)$ as
\begin{equation}
\mathbf{s}_\zeta(x_k) =\left\{
\begin{split}
&-\frac{\nabla f(x_k)}{|\nabla f(x_k)|}
-\zeta\frac{\nabla \Phi(x_k)}{|\nabla \Phi(x_k)|},  ~~~~~&\text{if}~~\nabla \Phi(x_k) \neq 0;\\
& -\nabla f(x_k) ,                              &\text{if}~~\nabla \Phi(x_k) = 0.
\end{split}\right.
\label{eq: vector_s_multiple}
\end{equation}
\subsection{Main contributions}
\label{sec:results}
Theory, computational algorithms, and applications are the three pillars of an optimization method. Just as the title suggests, this manuscript mainly focuses on the algorithms and applications of GDAM. We also give essential theoretical guarantees on the global convergence of the method for a complete presentation. 

\vskip 2mm
We summarize the main contributions of this work:
\begin{itemize}
    \item [1.] We present \eqref{eq:GDAM}, which is derived from our previous work \cite{Chen}, whereby the equivalence of the two methods is not obvious.
    \item [2.] We give fundamental theoretical guarantees to the method using a continuous-time dynamical systems approach. 
    \item [3.] We present two computational algorithms based on \eqref{eq:GDAM}, a vanilla and an accelerated implementation.
    \item [4.] We apply the developed algorithms to two engineering applications, node-based shape optimization and sensor network localization, which demonstrate the practical usefulness of GDAM in solving large-scale and difficult optimization problems.
\end{itemize}

\vskip 2mm

\subsection{Organization of paper}
In section \ref{sec:related}, we first review related works focusing on algorithm and application aspects. In section \ref{sec:derivation}, we derive the method for single inequality constrained problems. Illustrative examples are presented in section \ref{sec:preliminary} for an informal but intuitive presentation of the method. We give fundamental theoretical guarantees of the method in \ref{sec: global_behavior} and show its generalization to multiple constraints in section \ref{sec:multiple_constraint}. We present computational algorithms of GDAM in section \ref{sec:algorithms} and show preliminary numerical experiments. Two engineering applications, shape optimization and sensor network localization, are presented in sections \ref{sec:shapopt} and \ref{sec:snl}. We conclude the work with a discussion in section \ref{sec:conclusion}.


\section{Related works and applications}
\label{sec:related}

\subsection{Gradient descent method}
The gradient descent method, originally proposed by Cauchy, is a first-order method for unconstrained optimization that uses the negative gradient of the objective function for the variable update  \cite{lemarechal2012cauchy}. In a canonical form, it writes
\begin{equation}
    x_{k+1} = x_k - \alpha_k \nabla f(x_k),
\end{equation}
The \textit{direction} of the gradient descent is the steepest descent direction in Euclidean norm. To see this, we write the first-order Taylor expansion at the current iterate $x_k$ of the objective function,
\begin{equation*}
f(x_k + d) \approx f(x_k) + \langle \nabla f(x_k), d\rangle.
\end{equation*}
\textcolor[rgb]{0,0,0}{where ``$\langle\,,\rangle$" denotes the inner product in Euclidean space.}
The steepest descent direction $d$ is found by the optimization problem
\begin{equation*}
\begin{split}
&\textnormal{minimize}~~~ \langle \nabla f(x_k), d\rangle \\
&\textnormal{subject to}~~~ | d| = 1.
\end{split}
\end{equation*}
From the Cauchy-Schwarz inequality, we obtain
\begin{equation}
 d(x_k) = - \frac{\nabla f(x_k)}{|\nabla f(x_k)|}, 
 \label{eq:steepest_descent}
\end{equation}
which is the negative and normalized objective gradient.
Comparing (\ref{eq:steepest_descent}) and (\ref{eq:GDAM}),
\begin{equation*}
\mathbf{s}_\zeta (x_k) =-\frac{\nabla f(x_k)}{|\nabla f(x_k)|}
-\zeta\frac{\nabla g(x_k)}{|\nabla g(x_k)|}, ~~ \zeta \in [0,1),
\end{equation*}
which linearly combines the negative and normalized objective and constraint gradient, we consider the present method: a gradient descent akin method for inequality constrained optimization problems. The present method can be explained intuitively: For an inequality constrained problem, we use the objective gradient descent direction to minimize $f(x)$, and use the constraint gradient descent direction to minimize $g(x)$ (so as to maintain feasibility by moving away from the constraint boundary). Since the gradient directions are normalized to $1$, a parameter $ 0 \leq \zeta < 1$ ensures that the computed direction is a descent direction of the objective function. 

\subsection{Interior-point methods}

Interior-point methods (IPMs), which are mostly based on the Newton method, are among the most competitive methods for constrained optimization problems. The signature of IPM is the existence of continuously parameterized families of approximate solutions that converge to the exact solution asymptotically \cite{Forsgren}. IPMs find a wide variety of applications of convex and nonconvex optimizations in broad fields. There are a vast amount of excellent works that have been devoted to IPM. A comprehensive review of this class of methods is certainly beyond the scope of the present paper, however we refer the interested reader to \cite{Forsgren}\cite{Potra}\cite{ye2011interior}, more recently, in \cite{gondzio2012interior}, and many other excellent optimization books. 

The present method results in trajectories that asymptotically converge to the central path for a particular IPM, the logarithmic barrier method \cite{fiacco1990nonlinear}. We introduce its connections/differences to the present method in the next. Consider a convex optimization problem of the form (\ref{eq:Optimization_Problem}), we start with its approximated unconstrained subproblem using the logarithmic barrier function $\Phi(x)$,
\begin{equation}
\textnormal{minimize} ~~~~  f(x) + \eta \Phi(x),
\label{eq:logarithmic_barrier}
\end{equation}
where the barrier parameter $\eta$ is a positive parameter. As $\eta \rightarrow 0$, the solution of the approximated problem converges to the original one. The central path is characterized by the set of points that satisfy the necessary and sufficient conditions \cite{boyd}: 
\begin{equation}
0 = \nabla f(x^\star) + \eta \nabla \Phi(x^\star), ~~x^\star \in \Omega_-,
\label{eq:barrier_central_point}
\end{equation}
where $\Omega_- = \{x: g_i(x) < 0, i = 1,...,m\}$. The conditions (\ref{eq:barrier_central_point}) are interpreted as a modified KKT system in the literature \cite{boyd}\cite{byrd2000trust}\cite{Forsgren}. The barrier method finds an approximated solution for the original problem by 1) iteratively decreasing the barrier parameter $\eta$, and 2) in each iteration, solving the subproblem (the modified KKT system) defined by (\ref{eq:barrier_central_point}) using the Newton method. 
Therefore, the barrier parameter $\eta$ can be considered a central path parameter.

In the present method, we do not parameterize the central path. Instead, we normalize the objective and constraint gradients. For optimization problems with a single inequality constraint, we propose the \textit{normalized central path condition} as
\begin{equation}
\frac{\nabla f(x)}{|\nabla f(x)|} + \frac{\nabla g(x)}{|\nabla g(x)|}  = 0.
\label{eq:normalized_central_path_condition_single_constraint}
\end{equation}

The condition (\ref{eq:normalized_central_path_condition_single_constraint}) is used in section \ref{sec: global_behavior} for the analysis. A generalization of the \textit{normalized central path condition} from single inequality to multiple inequalities is introduced by the use of the logarithmic barrier function (see section \ref{sec:multiple_constraint}) as
\begin{equation}
\frac{\nabla f(x)}{|\nabla f(x)|} + \frac{\nabla \Phi(x)}{|\nabla \Phi(x)|}  = 0.
\label{eq:normalized_central_path_condition}
\end{equation}
Compared with the barrier method, the path parameter $\eta$ has vanished. The condition (\ref{eq:normalized_central_path_condition}) characterizes the central path, which differs from (\ref{eq:barrier_central_point}), which instead characterizes a point on the central path.

To see the connections between GDAM and the interior-point method, we first write the gradient descent direction $d_{\Phi,\eta}$ for a barrier $\eta$-subproblem \eqref{eq:logarithmic_barrier},
\begin{equation}
    d_{\Phi,\eta} = - \nabla f(x) - \eta \Phi(x).
\label{eq:gd_barrier}    
\end{equation}
A straightforward way to design a first-order method is to follow the double-loop structure of the second-order IPM: 1) iteratively decrease the barrier parameter $\eta$ in the outer loop, and 2) solve the subproblem \eqref{eq:logarithmic_barrier} by the gradient descent \eqref{eq:gd_barrier} in the inner loop. Unfortunately, this approach has long been known to be computationally impractical, see, e.g., the discussions in \cite[p. 469]{luenberger2021linear}. Roughly speaking, when $\eta$ is very small, $d_{\Phi,\eta}$ tends to result in an optimization trajectory that travels alongside the boundary of the constraints and suffers severely from poor conditioning. 

Next, we scale the GDAM search direction $s_\zeta$ with $|\nabla f(x_k) |$,
\begin{equation}
|\nabla f(x_k) | s_\zeta = -\nabla f(x_k) - \zeta \frac{|\nabla f(x_k) |}{|\nabla \Phi(x_k) |} \nabla \Phi(x_k ).
\label{eq:gdam_barrier_connection}
\end{equation}
Comparing \eqref{eq:gdam_barrier_connection} with \eqref{eq:gd_barrier}, it is then clear that the (scaled) GDAM suggests a \textit{dynamic} computation of the barrier parameter at \textit{each step} $k$,
\begin{equation}
 \eta(x_k) =\zeta \frac{|\nabla f(x_k) |}{|\nabla \Phi(x_k) |}.
 \label{eq:barrier_zeta}
\end{equation}
GDAM results in an optimization trajectory that travels alongside the central path (within a neighborhood relative to the parameter $\zeta$) \textemdash a typical behavior of a path-following method. Therefore, we consider GDAM as a first-order interior-point method. As one of the main contributions of this work, we will show that GDAM, in contrast to the sequential application of gradient descent to the barrier subproblems, is a computationally practical optimization method.

\subsection{Dynamical systems approaches} 
Dynamical systems approaches have been used to study optimization methods in many works, and our literature review could not be exhaustive. Extensive studies on the connections between interior-point flows with linear programming methods can be found in \cite[Chapter~4]{helmke2012optimization} and the references therein. Nonlinear dissipative dynamical systems are studied in \cite{attouch2000heavy} in view of unconstrained optimization. \cite{su2014differential} studies the celebrated Nesterov's accelerated gradient method using a dynamical system as the analysis tool. Recently, dynamical systems have been used to study optimization algorithms for solving problems arising from machine learning applications, see, e.g., \cite{arora2018optimization}\cite{chizat2018global}\cite{jordan2018dynamical}\cite{wilson2021lyapunov}. In some literature, optimization methods that use dynamical systems are called trajectory methods. These methods construct optimization paths in a way so that one or all solutions to the optimization problem are \textit{a priori} known to lie on these paths \cite{diener1995trajectory}. Typically, these optimization paths are solution trajectories to ODE of first or second-order. Trajectory methods are mainly studied for unconstrained optimizations for finding local solutions \cite{behrman1998efficient}\cite{botsaris1978differential}, and global solutions \cite{griewank1981generalized}\cite{snyman1987multi}. Studies for constrained optimization are, however, very limited, see \cite{ali2018trajectory}\cite{boct2020primal}\cite{wang2003unified} and the references therein. {In this work, we give basic theoretical guarantees of the present method using a dynamical system's perspective \textemdash we show that the method is globally convergent to first-order stationarities. 

\subsection{Feasible direction methods}
The method of feasible directions (MFD) dates back to the 1960's by the work of Zoutendijk \cite{zoutendijk1960methods} and has enjoyed fruitful developments for decades. MFDs have been especially popular in the engineering community because of the importance of ending up with a design that satisfies the hard specifications expressed by a set of inequalities \cite{chen2000methods}. The general idea behind the MFD is to move from one feasible design to an improved feasible design iteratively so that a local solution can be found \cite{Arora}. In the present method, the search direction must not be a feasible direction. The idea behind the method design is to find a search direction that approaches the central path while maintaining a descent direction of the objective function. Due to this major difference, we do not categorize our method as a method of feasible directions.

In the present method, we compute a search direction that uses normalized gradients. In \cite{stander1993new}, the authors also present a feasible direction method that applies normalized gradients. Their work is then continued and further developed in \cite{de1994feasible} and \cite{stander1995robustness}. In these works, the active-set strategy is used. The common idea is to formulate a linear system under given input criteria on a chosen working set of (active) constraints, and a feasible descent search direction is obtained by solving the linear system. In the present method, we use a barrier function based formulation to treat multiple inequality constraints. The search direction is computed as a linear combination of the negative and normalized objective and barrier gradient.

\subsection{Shape optimization}
As a subset of design optimization, shape optimization is characterized by a very large or even infinite number of design variables that describe the varying boundary in the optimization process. Introductions to shape optimization are given in \cite{haslinger2003introduction}\cite{sokolowski1992introduction}. Shape optimization is distinct from another well-known problem in design optimization: topology optimization \cite{bendsoe2013topology} (sometimes referred to as the homogenization method \cite{allaire2012shape}). The main difference is that the topology optimization method removes smoothness and topological constraints in shape optimization, which results in different optimization formulations. Many topology optimization problems can be formulated in an (equivalent) convex optimization problem, while shape optimizations are typically nonlinear, and
often nonconvex \cite{hoppe2007adaptive}. This difference partially contributes to the fact that there are successful implementations of IPM for large-scale\footnote{In this work, we refer large-scale optimizations to problems that have a high-dimensional variable space.} topology optimization problems \cite{jarre1998optimal}\cite{kennedy2015large}\cite{kocvara2016primal}\cite{maar2000interior}, but only a few works have presented a shape optimization that uses an IPM as the optimizer \cite{antil2007path}\cite{herskovits2000shape}. In the latter works, the size of the shape optimization problem is only moderate so that the power of IPM is not fully exploited. One of the most successful methods for nonlinear topology optimization is the method of moving asymptotes (MMA) that was introduced by Svanberg in 1987 \cite{svanberg1987method}. In each iteration, MMA generates and solves an approximated convex problem related to the original one. For shape optimization, however, there is as yet no literature that discusses a large-scale problem using MMA. 

Shape optimization is a subject frequently considered in multidisciplinary design optimization (MDO) \cite{haftka1992options}. To design complex engineering systems, MDO considers multiple disciplines and their interactions. For shape optimization, ongoing efforts are devoted to the computation of the shape gradient for coupled disciplines, for example, for steady-state coupled problems \cite{kenway2014scalable}\cite{najian2020partitioned} and for transient coupled problems \cite{korelc2009automation}. While the development for the coupled-gradient is challenging, the reward is accurate derivatives and massive reductions in computational cost, which are essential for large-scale optimizations \cite{martins2013multidisciplinary}. In many practical circumstances, the shape geometry is reparameterized with finitely many parameters. A shape reparameterization is typically needed if there are producibility restrictions \cite{liedmann2020shape} or if the initial geometric variable is not differentiable \cite{hwang2015modular}. Well-parameterized shape geometry often allows the application of standard optimization approaches \cite{frohlich2019geometric}\cite{hicken2010aerodynamic}\cite{pyopt-paper} or Newton-Krylov type methods \cite{dener2015comparison}. On the other hand, the achievable shape is limited and dependent on the chosen parameterization. High-fidelity shape optimizations, which are directly based on finite element meshes \cite{linkwitz1971einige}\cite{zienkiewicz1977finite} or level-set methods \cite{allaire2004structural}\cite{sethian1999level}, exploit the largest shape space possible for real-life problems but lead to challenging optimization problems \cite{fepponnull}\cite{luft2020efficient}.

Manually deriving and implementing shape derivatives can be a laborious and error-prone task. A prominent approach to tackle this challenge is Algorithmic Differentiation (AD) \cite{griewank2008evaluating}, which computes derivatives of a function given as a computer program. In the context of shape optimization, AD has been successfully implemented in open-source solvers, such as OpenFOAM\cite{towara2013discrete}, SU2\cite{albring2016efficient}, and FEniCS\cite{dokken2020automatic}, enabling shape design optimization practice on increasingly complex problems. Another ongoing research is the computation of shape Hessians, which are complex objects even for moderate problems. Recently, several works compute approximated shape Hessians and use a Newton-based method for the design optimization \cite{schillings2011efficient}\cite{schmidt2013three}. In general, large-scale shape optimization is mainly performed using gradient descent type methods so far \cite{schulz2015}. In engineering practice, a large number of constraints may be considered. The lack of literature in this regard has motivated our development of MSDM in \cite{Chen}. In the present work, we simplify its computational framework to a gradient descent akin method. As an important result, the implementation effort and computational cost are reduced significantly. It opens the possibility of shape optimization to a wider range of applications.

\subsection{Sensor network localization via semidefinite programming relaxation}

Wireless sensor networks (WSNs), which consist of low-cost, low-energy, and multi-functional sensors that can communicate over short distances, provide many opportunities for monitoring and controlling the physical environment when the sensors are wireless linked and deployed in large numbers \cite{akyildiz2002wireless}. WSNs find a wide range of applications in environmental, industrial, urban, health, and other sectors, and we refer to comprehensive surveys in \cite{kandris2020applications}\cite{puccinelli2005wireless}. For many applications, awareness of the sensor locations is crucial for a meaningful interpretation of the gathered sensing data, see, e.g.,  \cite{buehrer2018collaborative}\cite{leonard2007collective}\cite{raty2010survey}\cite{sun2005reliable}. Sensor network localization (SNL), which estimates sensor locations based on pairwise distance, is therefore considered one of the key enabling technologies for WSN applications. 

SNL is a challenging optimization problem, as it is a nonconvex problem and is generally intractable to find the global solution. Among various solution methods, convex relaxation based on semidefinite programming (originally introduced in Biswas-Ye \cite{biswas2004semidefinite}) is regarded as one of the most prominent approaches for global localization \cite{chowdhury2016advances}. In theory, the semidefinite programming (SDP) relaxation provides the exact solution to an SNL problem if the problem is uniquely localizable \cite{so2007theory}. The challenge that has hindered the practical application of SDP relaxation to large-scale SNL problems is primarily computational\footnote{Indeed, scalability is a major challenge to the successful application of semidefinite programming to many large-size practical problems.}. Well-established SDP solvers that are based on second-order interior-point methods (\cite{andersen2000mosek}\cite{benson2008algorithm}\cite{polik2007sedumi}\cite{toh2012implementation}\cite{yamashita2012latest}) can solve SNL problems with up to a few hundreds of sensors to arbitrary high accuracy, but they are unable to solve larger problems efficiently. This has motivated the development of advanced modelling approaches to alleviate the computational difficulty, see, e.g., \cite{biswas2006distributed}\cite{kim2009exploiting}\cite{wang2008further}. Meanwhile, in the last two decades, a number of scalable SDP solvers were successfully developed. For example, the packages PENNON \cite{kovcvara2003pennon} and SDPNAL+ \cite{sun2020sdpnal+} are based on the Augmented Lagrangian method, and SCS \cite{odonoghue21} is based on the Alternating Direction of Multipliers method (ADMM).
A comprehensive review of the latest advancements in the scalability of SDP solvers is given in \cite{majumdar2020recent}. In this work, we apply GDAM to solve the Biswas-Ye SDP relaxation of SNL problems. We use the classical logarithmic barrier function approach for the semidefinite cone constraint, and thus our method to solve SDP problems is general. Our computational experiments show that GDAM, when practically implemented, is very competitive in finding moderate accurate solutions to SNL problems of moderate size and is capable of solving very large problems (over 5000 sensors) that are intractable for comparing solvers.

\section{Deriving the search direction for single inequality constrained optimizations}
\label{sec:derivation}
In this section, we show the consistent derivation of \eqref{eq:GDAM} from our previous work \cite{Chen} considering a single inequality constraint. First, we review the basic ideas of the modified search direction method and then show the derivation. For the remainder of the paper, we use $\nabla f$ and $\nabla g$ instead of $\nabla f(x)$ and $\nabla g(x)$ to lighten the notation when it is clear from the context.

\textcolor[rgb]{0,0,0}{In MSDM, at each iterate $x$, we construct a sensitivity matrix that contains the normalized objective and constraint gradient,
\begin{equation}
\mathbf{m}(x) = \begin{bmatrix}  \frac{\nabla f^T}{|\nabla f|} \\ \frac{\nabla g^T}{|\nabla g|} \end{bmatrix}.
\label{eq: sensitivity_matrix}
\end{equation}}
\textcolor[rgb]{0,0,0}{The total differential of the objective and constraint function at $x$ are
\begin{equation}
\begin{split}
df = \nabla f^T dx, \\
dg = \nabla g^T dx.
\end{split}
\end{equation}}
A perspective from the input-output system established by \textcolor[rgb]{0,0,0}{the matrix $\mathbf{m}(x)$ gives
\begin{equation}
\begin{pmatrix}  \frac{df}{|\nabla f|} \\\frac{dg}{|\nabla g|} \end{pmatrix} = \mathbf{m}(x) d x.
\label{eq:input_output_smatrix}
\end{equation}}
Applying SVD to \textcolor[rgb]{0,0,0}{the matrix $\mathbf{m}(x)$, we obtain
\begin{equation}
\mathbf{m}(x) = \mathbf{U} \mathbf{\Sigma} \mathbf{V}^T = \sum_{i = 1}^{min(2,n)} \mathbf{\sigma}_{i} \mathbf{u}_i \mathbf{v}_i^T.
\label{eq: svd_m}
\end{equation}}
Thus, an orthonormal bases set $\mathbf{v}_i, ~i = 1,2, \mathbf{v}_i \in \mathbb{R}^n$ is obtained. Each $\mathbf{v}_i$ can be used as a base search direction for the \textcolor[rgb]{0,0,0}{variable update}, 
and $\mathbf{v}_1$ and $\mathbf{v}_2$ are defined as follows:
\begin{itemize}
\item[-] $\mathbf{v}_1$: by taking $\delta \mathbf{v}_{1}$ as the \textcolor[rgb]{0,0,0}{variable update}, we obtain a change in objective as well as in constraint function $[\frac{df}{|\nabla f|}, \frac{dg}{|\nabla g|}]^T = \sigma_{1} \delta \mathbf{u}_{1} $, which is a \textbf{decrease} in the objective function and an \textbf{increase} in the constraint function.

\item[-] $\mathbf{v}_2$: by taking $\delta \mathbf{v}_{2}$ as the \textcolor[rgb]{0,0,0}{variable update}, we obtain a change in the objective as well as in the constraint function $[\frac{df}{|\nabla f|}, \frac{dg}{|\nabla g|}]^T = \sigma_{2} \delta \mathbf{u}_{2} $, which is a \textbf{decrease} in the objective function and a \textbf{decrease} in the constraint function. 	
\end{itemize}
It is worth mentioning that the \textcolor[rgb]{0,0,0}{update step} $\delta \mathbf{v}_1$ provides a similar result as the filter approach presented in \cite{fletcher2002nonlinear}, which tries to minimize the so-called bi-objective optimization problem with two goals of minimizing the objective function $f$ and the constraint violation $|g|$. 

With $\mathbf{v}_1$ and $\mathbf{v}_2$, we can then rewrite the normalized steepest descent direction $-\frac{\nabla f}{|\nabla f|}$ as
\begin{equation}
-\frac{\nabla f}{|\nabla f|} = \cos \alpha_1 \mathbf{v}_{1} + \cos \alpha_2 \mathbf{v}_{2},
\label{eq:rewrite_steepest_descent}
\end{equation}
where $\alpha_1$ is the angle between $\mathbf{v}_{1}$ and $-\frac{\nabla f}{|\nabla f|}$, and $\alpha_2$ is the angle between $\mathbf{v}_{2}$ and $-\frac{\nabla f}{|\nabla f|}$. The modified search direction reads
\begin{equation}
\mathbf{s}_c = \cos \alpha_1 \mathbf{v}_{1} + c \cdot \cos \alpha_2 \mathbf{v}_{2},
\label{eq:ds_scop_svd}
\end{equation}
where $c \geq 1$ is introduced to enlarge the contribution of $\mathbf{v}_2$ for the variable update.

In the following, we show the derivation of \eqref{eq:GDAM} from \eqref{eq:ds_scop_svd}. By \eqref{eq: svd_m}, we have
\begin{equation}
\mathbf{m} \mathbf{m}^T = \mathbf{U} \mathbf{\Sigma} \mathbf{\Sigma}^T \mathbf{U}^T,
\end{equation}
which is a diagonalization of the symmetric matrix $\mathbf{m} \mathbf{m}^T$. We have
\begin{equation}
\mathbf{m} \mathbf{m}^T= \begin{pmatrix}  \frac{\nabla f^T}{|\nabla f|} \\ \frac{\nabla g^T}{|\nabla g|} \end{pmatrix} 
\begin{pmatrix}
\frac{\nabla f}{|\nabla f|} & \frac{\nabla g}{|\nabla g|}
\end{pmatrix} = 
\begin{pmatrix}
1 & \cos \theta \\
\cos \theta & 1
\end{pmatrix},
\end{equation}
where $\theta$ is the angle between $\nabla f$ and $\nabla g$. The eigenvalues of $\mathbf{m} \mathbf{m}^T$ are
\begin{equation}
\begin{split}
\lambda_1 = 1 - \cos \theta, \\
\lambda_2 = 1 + \cos \theta. \\
\end{split}
\end{equation}
By the definition of $\mathbf{v}_1$ and $\mathbf{v}_2$ in MSDM, the eigenvectors of $\mathbf{mm}^T$ can be factorized,
\begin{equation}
\begin{split}
\mathbf{u}_1 &= \begin{pmatrix}
-\frac{\sqrt{2}}{2}, & \frac{\sqrt{2}}{2}
\end{pmatrix}^T, \\
\mathbf{u}_2 &= \begin{pmatrix}
-\frac{\sqrt{2}}{2}, & -\frac{\sqrt{2}}{2}
\end{pmatrix}^T. \\
\end{split}
\end{equation}
The respective singular values are
\begin{equation}
\begin{split}
\sigma_1 &= \sqrt{1- \cos \theta}, \\
\sigma_2 &= \sqrt{1+ \cos \theta}. \\
\end{split}
\end{equation}
By the definition of SVD,
\begin{equation}
\begin{pmatrix}  \frac{\nabla f^T}{|\nabla f|} \\ \frac{\nabla g^T}{|\nabla g|} \end{pmatrix} = 
\begin{pmatrix}
-\frac{\sqrt{2}}{2} & -\frac{\sqrt{2}}{2} \\
\frac{\sqrt{2}}{2} & -\frac{\sqrt{2}}{2}
\end{pmatrix}
\begin{pmatrix}
\sqrt{1- \cos \theta} & 0 \\
0 & \sqrt{1+ \cos \theta} 
\end{pmatrix}
\begin{pmatrix}
\mathbf{v}_1^T \\
\mathbf{v}_2^T
\end{pmatrix}.
\end{equation}

Therefore, we obtain
\begin{equation}
\begin{split}
\mathbf{v}_1&=\frac{1}{\sqrt{2-2\cos\theta}}\left( -\frac{\nabla f}{|\nabla f|}+\frac{\nabla g}{|\nabla g|}\right),
\\
\mathbf{v}_2&=\frac{1}{\sqrt{2+2\cos\theta}} \left(-\frac{\nabla f}{|\nabla f|}-\frac{\nabla g}{|\nabla g|}\right),
\end{split}
\label{eq:base_vectors}
\end{equation}
where $\theta$ is the angle between the objective function gradient $\nabla f$ and the constraint function gradient $\nabla g$ \textcolor[rgb]{0,0,0}{(full details in Appendix A)}. With $\theta$ we also have
\begin{equation}
\begin{split}
\cos \alpha_1=-\langle \frac{\nabla f}{|\nabla f|},\mathbf{v}_1\rangle=\frac{\sqrt{1-\cos\theta}}{\sqrt{2}},\\ \cos \alpha_2=-\langle\frac{\nabla f}{|\nabla f|},\mathbf{v}_2\rangle=\frac{\sqrt{1+\cos\theta}}{\sqrt{2}}.
\end{split}
\label{eq:cos_alphas}
\end{equation}
Inserting (\ref{eq:base_vectors}), (\ref{eq:cos_alphas}) into (\ref{eq:ds_scop_svd}) we have
\begin{equation}
\mathbf{s}_c= -\frac{\nabla f}{|\nabla f|}-\frac{(c-1)}{2} \left(\frac{\nabla f}{|\nabla f|}+\frac{\nabla g}{|\nabla g|}\right).
\end{equation}
As we are mainly interested in the direction of the vector field $\mathbf{s}_c$, we can rewrite it as
\begin{equation}
\mathbf{s}_\zeta =-\frac{\nabla f}{|\nabla f|}
-\zeta\frac{\nabla g}{|\nabla g|},
\label{vector_s}
\end{equation}
with $\zeta = \frac{c-1}{c+1}$. With $c \in [1,+\infty)$ we have $\zeta \in [0, 1)$ and thus we get the present search direction \eqref{eq:GDAM}.

\section{Illustrative examples and intuition}
\label{sec:preliminary}
In this section, we demonstrate GDAM with two illustrative examples and provide intuitive explanations for its optimization behavior. 

\subsection{The behavior of the method with a fixed $\zeta$}
\label{sec:nonconvex_example}
We consider a quadratically constrained nonconvex optimization problem,
\begin{equation}
\begin{split}
&\textnormal{minimize}~~~	f(x_1,x_2)= 1.5x_1^2 + x_2^2 - 2x_1 x_2 + 2x_1^3 + 0.5 x_1^4, \\
&\textnormal{subject to}~~  g(x_1,x_2)= x_1^2 - x_2 -2.2 \leq 0. \\
\end{split}
\label{eq:2d_nonconvex_obj_quadratic_constraint}
\end{equation}

The problem has multiple local solutions due to nonconvexity. Figure \ref{fig:2d_nonconvex_obj_quadratic_constraint} shows that there are two local solutions to the considered problem in the depicted variable space:
\begin{itemize}
    \item [1.] A critical point of $f(x)$, illustrated as a star, lies inside the feasible set.
    \item [2.] A KKT solution, illustrated as a dot, is located at the boundary of the constraint.
\end{itemize}

It can be observed that the optimization trajectories all converge to the critical points (shown as stars). Moreover, the intersection points of the two left trajectories and the constraint boundary are approximate KKT points. While the former case is a common behavior of gradient flow approaches, the latter case is unique to the present method.

\begin{figure}[h]
	\centering
	\begin{footnotesize}
		\includegraphics[width=12.5cm]{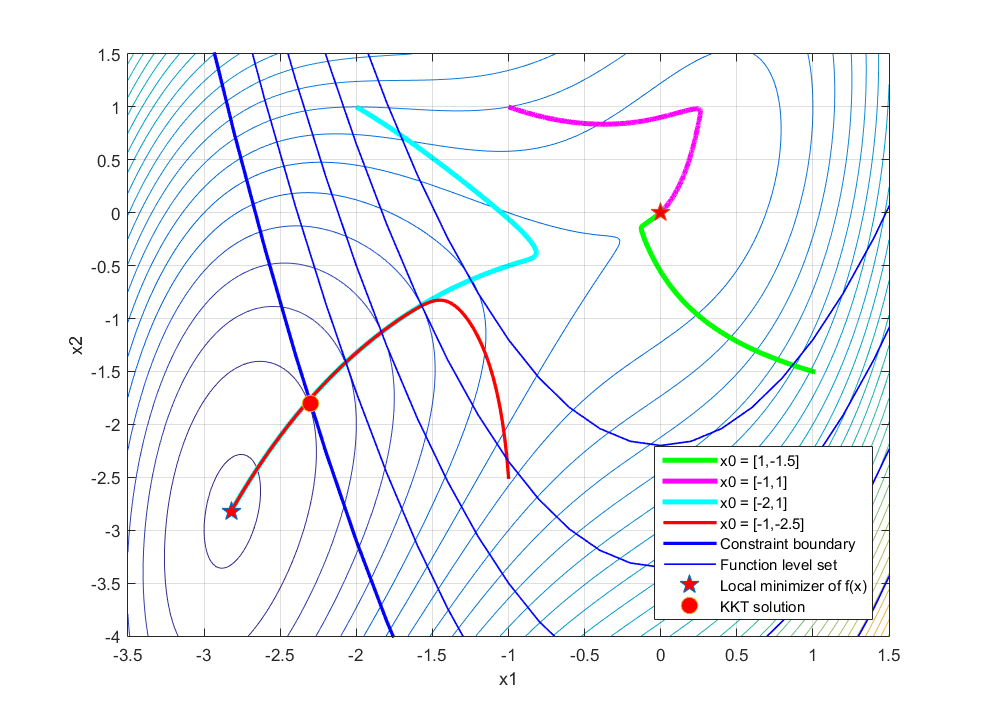}
		\caption{Optimization trajectories for problem \eqref{eq:2d_nonconvex_obj_quadratic_constraint} with $\zeta = 0.95$ and different initializations. Depending on the initialization, different solutions to the optimization problem are found.}
		\label{fig:2d_nonconvex_obj_quadratic_constraint}
	\end{footnotesize}
\end{figure}

\subsection{The behavior of the method in terms of $\zeta$}
\label{sec:2d_analytic}
We show the behavior of the method in terms of different $\zeta$ by solving a 2D optimization problem analytically.
The optimization problem reads,
\begin{equation}
\begin{split}
&\textnormal{minimize}~~~	f(x_1,x_2)=\frac{1}{2}(x_1^2 + x_2^2), \\
&\textnormal{subject to}~~  g(x_1,x_2)= - x_2 + 10 \leq 0. \\
\end{split}
\label{eq:example_proof}
\end{equation}

\begin{figure}[h]
	\centering
	\begin{footnotesize}%
	\executeiffilenewer{fig/2d_linear_constraint_analysis_zeta.svg}{fig/2d_linear_constraint_analysis_zeta.pdf}%
	{inkscape -z -C --file=fig/2d_linear_constraint_analysis_zeta.svg %
		--export-pdf=fig/2d_linear_constraint_analysis_zeta.pdf --export-latex}%
	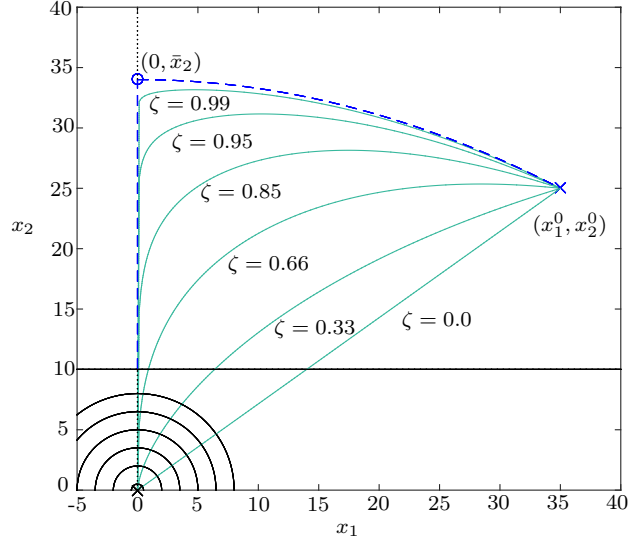%

		\caption{Optimization trajectories for the 2D linear constrained optimization problem \eqref{eq:example_proof} with different $\zeta$. The black circles show the objective function contours. The black line shows the constraint function. The dotted line is the central path. As $\zeta \rightarrow 1^{-}$, a portion of the trajectory $\Gamma^\zeta$ based on the present search direction field (\ref{vector_s}) converges to the central path.}
		\label{fig:2d_linear_constraint_analysis}
	\end{footnotesize}
\end{figure}

The present search direction field $\mathbf{s}_\zeta$ reads
\begin{equation}
\mathbf{s}_\zeta = \frac{-1}{\sqrt{x_1^2 + x_2^2}}\left\lbrace  x_1,  x_2 - \zeta \sqrt{x_1^2 + x_2^2} \right\rbrace^T.
\end{equation}
We define an initialization as $(x_1^0, x_2^0)$. Let $\bar{x}_2 = \frac{1}{2} \left(x_2^0 + \sqrt{(x_1^0)^2 + (x_2^0)^2}\right)$, then, the trajectory $\Gamma^\zeta$ of the present search direction field $\mathbf{s}_\zeta$ is 
\begin{equation}
\tag{$\Gamma^\zeta$}
x_2 + \sqrt{x_1^2 + x_2^2} = 2 \bar{x}_2 \left|\frac{x_1}{x_1^0}\right|^{1- \zeta}.
\label{eq:Gamma_zeta}
\end{equation}
Let  $(x_{1,\zeta}, x_{2,\zeta})$ be a point on the trajectory $\Gamma^\zeta$ with a maximal $x_2$ component, then we have
\begin{equation}
(x_{2,\zeta})^\zeta=\frac{2 \zeta}{1+\zeta}\frac{\bar{x}_2}{|x_1^0|^{1-\zeta}}
\left(\frac{\sqrt{1-\zeta^2}}{\zeta}\right)^{1-\zeta},
\end{equation}
and
\begin{equation}
|x_{1,\zeta}|=\frac{\bar{x}_2}{\zeta}\sqrt{1-\zeta^2}. 
\end{equation}
Let $\zeta \rightarrow 1^{-} $, then, $x_{1,\zeta} \rightarrow 0$, $x_{2,\zeta} \rightarrow \bar{x}_2$, the trajectory $\Gamma^\zeta$ will converge to the curve
$\Gamma$ that is a union of the parabola
\begin{equation}
x_1^2 =  4 \bar{x}_2^2 - 4 x_2 \bar{x}_2, ~~~~ x_1\in (0, x_1^0) (~~~\mbox{or} (x_1^0,0),)
\end{equation}
and the interval $(0,\bar{x}_2)$ on $x_2$-axis as is shown in figure \ref{fig:2d_linear_constraint_analysis} (full derivation details in Appendix \ref{appendix:1}).

We now give an estimation of the intersection point $x_\zeta^\sharp$ of the trajectory \eqref{eq:Gamma_zeta} with the constraint boundary,
\begin{equation}
	\left\{
	\begin{split}
	&x_{\zeta,2}^\sharp + \sqrt{(x_{\zeta,1}^\sharp)^2 + (x_{\zeta,2}^\sharp)^2} = 2 \bar{x}_2 \left|\frac{x_{\zeta,1}^\sharp}{x_1^0}\right|^{1- \zeta},\\
	&x_{\zeta,2}^\sharp = 10;
	\end{split}
	\right.
\end{equation}
For a simple presentation, we first assume that $x_1^0 > 0$. By
$$\frac{dx_1}{dt} = \frac{-x_1}{\sqrt{x_1^2 + x_2^2}},$$
the optimization trajectory is monotonically decreasing as $x_1 >0$, and thus $x_{\zeta,1}^\sharp < x_{1,\zeta}$. Furthermore, $\frac{dx_1}{dt}|_{x_1=0} = 0$ ensures that $x_1(t) \geq 0$ holds for any $t$. Combining both arguments, we obtain an estimate for the intersection point $x_{\zeta}^\sharp$,
\begin{equation}
	\left\{
	\begin{split}
	&    0\leq x_{\zeta,1}^\sharp \leq \frac{\bar{x}_2}{\zeta}\sqrt{1-\zeta^2},\\
	&x_{\zeta,2}^\sharp = 10.
	\end{split}
	\right.
\end{equation}
Obviously, $x_1^0 < 0 $ allows a similar study. Thus we obtain 
\begin{equation}
    | x_{\zeta}^\sharp - x^\star| \leq \frac{\bar{x}_2}{\zeta}\sqrt{1-\zeta^2},
\label{eq:estimate_analytic_sol}    
\end{equation}
where $x^\star = (0,10)$ is the KKT solution. We note that the error bound \eqref{eq:estimate_analytic_sol} is derived specifically for the problem \eqref{eq:example_proof} and is a very rough estimation that is based on a special point $(x_{1,\zeta}, x_{2,\zeta})$. In the following section, we will give a general result regarding such an error bound, as well as other fundamental theoretical guarantees.

\section{Global convergence}
\label{sec: global_behavior}

In this section, we give fundamental theoretical guarantees of the method regarding global convergence. We shall first analyze the single constraint optimization problem (SCOP), whose results can then be easily generalized to multiple constraints via the logarithmic barrier function (see next section). 

For now, let's consider
\begin{equation}
\tag{SCOP}
\begin{split}
&\textnormal{minimize}~~~  f(x), \\
&\textnormal{subject to}~~ g(x) \leq 0. \\
\end{split}
\label{eq:single_constrained_optimization_problem}
\end{equation}
Recall the iterative update formula of \eqref{eq:GDAM},
\begin{equation*}
\tag{GDAM}
\begin{split}
    &x_{k+1} = x_k + \alpha_k \mathbf{s}_\zeta (x_k), \\
    &\mathbf{s}_\zeta (x_k) = -\frac{\nabla f(x_k)}{|\nabla f(x_k)|}
-\zeta\frac{\nabla g(x_k)}{|\nabla g(x_k)|}, ~~ \zeta \in [0,1).
\end{split}  
\end{equation*}
We make use of the following dynamical system to study the convergence of the method,
\begin{equation}
\tag{DS}
\frac{d x}{d t} =-\frac{\nabla f(x)}{|\nabla f(x)|}
-\zeta\frac{\nabla g(x)}{|\nabla g(x)|}, ~~ \zeta \in [0,1).
\label{eq:ds_scop}
\end{equation}
In addition to the theoretical convergence guarantee, we report results on the global behavior of the continuous-time trajectory of the method, which can be useful for the design of practical algorithms.

\subsection{Assumptions and basic results}

For the analysis of the system \eqref{eq:ds_scop}, we make the following assumptions:
\begin{itemize}
\item[(A1)]  Coercive condition for the objective function $f(x)$,

$$\lim_{x\rightarrow \infty}f(x)=+\infty;$$

\item[(A2)] $\nabla f(x) \neq 0$ in the feasible set $\Omega=\{x: g(x) \leq 0\}$;

\item[(A3)] $\nabla g(x) \neq 0$ in the feasible set $\Omega$;

\item[(A4)]  $f$ and $g$ are twice continuously differentiable functions.
\end{itemize}
\vskip 2mm

If the function $f(x)$ is \textit{coercive} and continuous, then it has a global minimizer. The assumptions (A2) and (A3) are introduced due to our continuous-time analysis, as we need system \eqref{eq:ds_scop} to be well-defined. We will discuss the case where (A2) is excluded in subsection \ref{sec:global_conv_critical_points}. In practice, (A3) may not always be satisfied. To escape the critical points of constraint functions, we suggest using the gradient descent search direction (see also section \ref{sec:multiple_constraint}). 

Let $x(t; \zeta,x_0)$ be the solution of the system:

\begin{equation}\label{system}
\left\{
\begin{split}
&\frac{dx}{dt}= \textbf{s}_\zeta(x),\\
& x|_{t=0} = x_0.
\end{split}\right.
\end{equation}
with $x_0$ as the initial design in the feasible set and $T_{\zeta,x_0}$ as the \textit{maximal existence interval} of $x(t; \zeta,x_0)$ in the set $E \subset \mathbb{R}^n$, $$E := \{x: |\nabla f(x)| \neq 0, |\nabla g(x)| \neq 0, x \in \mathbb{R}^n \}.$$ 
\textcolor[rgb]{0,0,0}
{We refer to \cite{perko2008differential} for a formal definition of the maximal interval of existence}.  In the following, we show some basic properties of the trajectory  $x(t; \zeta,x_0)$.
\vskip 2mm

\subsubsection{Time derivative of $f(x(t,\zeta, x_0))$ and $g(x(t,\zeta, x_0))$}
We first observe the changes of the objective function $f(x)$ and constraint function $g(x)$ along the solution trajectory $x(t,\zeta, x_0)$. This is done by taking the time derivatives of $f(x(t,\zeta, x_0))$ and $g(x(t,\zeta, x_0))$,
\begin{equation}
\begin{split}
&\frac{d}{dt} \left( f(x(t;\zeta, x_0)) \right) = \frac{df(x)}{dx} \cdot \frac{dx(t;\zeta, x_0)}{dt} \\
&= \nabla f(x)^T \left(- \frac{\nabla f(x)}{|\nabla f(x)|} - \zeta \frac{\nabla g(x)}{|\nabla g(x)|} \right) \\
&= - |\nabla f(x)| \left( \frac{\nabla f(x)^T \nabla f(x)}{|\nabla f(x)| |\nabla f(x)|} + \zeta \frac{\nabla f(x)^T \nabla g(x)}{|\nabla f(x)| |\nabla g(x)|}  \right) \\
& = -|\nabla f(x)|\left(1 + \zeta \cos \theta(x) \right),
\end{split}		
\label{eq:deform_f}
\end{equation}
similarly,
\begin{equation}
\frac{d}{dt} \left(g(x(t;\zeta, x_0))\right) = - |\nabla g(x)| ( \cos \theta(x) + \zeta ),
\label{eq:deform_g}
\end{equation}
where 
\begin{equation}
   \cos \theta (x) =  \frac{\nabla f(x)^T \nabla g(x)}{|\nabla f(x)| |\nabla g(x)|}.
   \label{eq:cos_theta}
\end{equation} 
Obviously, \eqref{eq:deform_f}, \eqref{eq:deform_g}, and \eqref{eq:cos_theta} hold under Assumptions $(A1)-(A4)$ and as $x \in \Omega$. In the following, we derive basic theoretical results based on the above three equations.

\subsubsection{Boundedness of the solution $x(t;\zeta, x_0)$} First, it is easy to observe from \eqref{eq:deform_f} that $f(x)$ decreases along $x(t;\zeta, x_0)$ with a fixed $\zeta \in [0,1)$. Furthermore, there exists a global lower bound for the objective function $f(x)$ by the $coercive$ Assumption (1). We can show that the solution trajectory $x(t;\zeta, x_0)$ always stays in a bounded set that is dependent on the initialization $x_0$. The result is formalized as follows.

First, we define the \textit{bounded set} in relation to a feasible initialization $x_0$.
\begin{definition} [Bounded set $B_{f(x_0)}$]
   Given a feasible initialization $x_0$,
$$x_0 \in \Omega = \{x: g(x) \leq 0 \}. $$
The bounded set $B_{f(x_0)}$ is defined by
\begin{equation}
B_{f(x_0)} := \{x: f(x) \leq f(x_0)\}.
\label{eq:bounded_set}
\end{equation} 
\end{definition}

\begin{lemma} [Boundedness]
Suppose that assumptions (A1) - (A4) hold, and let $\zeta \in [0,1)$, $x_0 \in \Omega$. Then the trajectory $x(t; \zeta,x_0)$ always stays within the bounded set $B_{f(x_0)}$ associated with the feasible initialization $x_0$.
\label{lemma_boundedness}
\end{lemma}
\begin{proof}
~~Based on the deformation of the objective function $f(x)$ along the trajectory $x(t; \zeta,x_0)$ shown in \eqref{eq:deform_f}, and with $\zeta \in [0,1)$, $\cos \theta(x) \in [-1,1]$, we have
\begin{equation}
1 + \zeta \cos \theta(x) > 0.
\end{equation}
Therefore,
\begin{equation}
\frac{df(x(t;\zeta, x_0))}{dt} \leq 0.
\end{equation}
Thus, the proof is complete.
\end{proof}

\begin{remark}
\label{remark:boundedness}
Lemma \ref{lemma_boundedness} is a basic result that implies that the integral of \eqref{eq:deform_f} is always bounded under our assumptions, i.e.,
\begin{equation}
\int_0^{T_{\zeta, x_0}}|\nabla f(x(t; \zeta,x_0))|(1+\zeta \cos \theta (x) )dt = f(x_0)-f(x(T_{\zeta, x_0}; \zeta,x_0)) < \infty.
\label{eq:integral_f_boundedness}
\end{equation}
\end{remark}

\subsubsection{Lipschitz continuity of the vector field $\mathbf{s}_\zeta$}

The intersection of the feasible set $\Omega$ and the bounded set $B_{f(x_0)} $ defines the \textit{bounded feasible set}.
\begin{definition} [Bounded feasible set  $\Omega_{f(x_0)}$]
Given a feasible initialization $x_0$,
$$x_0 \in \Omega = \{x: g(x) \leq 0 \}. $$
The bounded feasible set $\Omega_{f(x_0)}$ is defined by
\begin{equation}
\Omega_{f(x_0)} := \{x: f(x) \leq f(x_0), x \in \Omega \}.
\label{eq:bounded_feasible_set}
\end{equation}
\end{definition}

In the bounded feasible set $\Omega_{f(x_0)}$, we show Lipschitz continuity of the present vector field $\mathbf{s}_\zeta$.
\vskip 2mm
	\begin{lemma}
		Suppose that assumptions (A1)-(A4) hold and let the bounded feasible set $\Omega_{f(x_0)}$ be non-empty.
		 Then, the vector field $\mathbf{s}_\zeta(x), ~\zeta \in [0,1]$ is Lipschitz continuous for $x \in \Omega_{f(x_0)}$, i.e.,
		\begin{equation}
		| \mathbf{s}_\zeta (x) - \mathbf{s}_\zeta (y) | \leq L | x - y|,
		\end{equation}
		\label{lemma_lipschitz}
where $L$ is a positive constant.
	\end{lemma}
	\begin{proof}
	    ~~By assumption (A4), the second derivatives of $f$ and $g$ are bounded in $\Omega_{f(x_0)}$. Therefore, $\nabla f$ and $\nabla g$ are Lipschitz continuous in $\Omega_{x_0}$. By assumption (A2) and (A3), we have
		\begin{equation*}
		| \nabla f(x)| \geq a, ~~ | \nabla g(x)| \geq b, ~~\forall x \in \Omega_{x_0},
		\end{equation*}
		for some positive numbers $a,b$. Therefore, $\frac{\nabla f(x)}{|\nabla f(x)|}$ and $\frac{\nabla g(x)}{|\nabla g(x)|}$ are Lipschitz continuous. Thus, we have a Lipschitz continuity for $\mathbf{s}_\zeta(x) = -\frac{\nabla f(x)}{|\nabla f(x)|}
		-\zeta\frac{\nabla g(x)}{|\nabla g(x)|}$ with $x \in \Omega_{f(x_0)}$.	
	\end{proof}

\vskip 2mm

\subsubsection{A cosine measure for the central path neighborhood}
\label{subsec: cosine}

For the studies of inequality constrained optimizations, the central path is recognized as ``a fundamental mathematical object'' \cite{bayer1989nonlinear}. For example, its total curvature is used in the study of the existence of strongly polynomial algorithms for linear programs in the context of IPMs \cite{allamigeon2018log}. The equation \eqref{eq:cos_theta} inspired us to give a cosine measure for the central path and its neighborhood. 

Recall the normalized central path condition \eqref{eq:normalized_central_path_condition_single_constraint},
\begin{equation*}
\frac{\nabla f(x)}{|\nabla f(x)|} + \frac{\nabla g(x)}{|\nabla g(x)|}  = 0,
\end{equation*}
which can be formulated as the angular relation between $\nabla f(x)$ and $\nabla g(x)$
\begin{equation}
\cos \theta(x) = \frac{\langle \nabla f(x),  \nabla g(x) \rangle }{| \nabla f(x) | | \nabla g(x) |} =  -1.
\label{eq:cos_theta_1}
\end{equation}

\indent Under our assumptions (A1)-(A4), $\cos \theta(x)$ is a continuously differentiable  function of $x$ in $\Omega_{f(x_0)}$. Additionally, we note that $\cos \theta(x)$ reaches its own minimum when $x$ is at the central path. Combining both arguments allows us to conveniently define a neighborhood of the central path using the level set of $\cos \theta(x)$.

\begin{definition} [$\mu-$neighborhood]
The $\mu-$neighborhood of a central path is defined as
$$\Theta_\mu=\{x:\cos\theta(x) <-\mu, x \in \Omega_{f(x_0)} \},~~\mu \in [0,1]. $$
\label{def:mu-neighborhood}
\end{definition}
Intuitively, $\Theta_\mu$ is a cone-like neighborhood around the central path. Obviously,  $\Theta_\mu$ shrinks to the central path as $\mu \rightarrow 1^-.$ With the definition of the $\mu-$neighborhood, we can show the following behavior of the optimization trajectory related to the constraint function $g(x)$.

\begin{lemma}
\label{lemma:deform_g}
	Suppose that assumptions (A1)-(A4) hold, $\zeta \in [0,1)$ is fixed, and $x_0 \in \Omega$, then the constraint function $g(x)$ decreases along the trajectory $x(t; \zeta,x_0)$ out of the neighborhood $\Theta_{\zeta}$, and increases in $\Theta_{ \zeta}$, with $\Theta_{ \zeta}$ defined as 
 \begin{equation}
     \Theta_\zeta=\{x:\cos\theta(x) <-\zeta, x \in \Omega_{f(x_0)} \}. 
 \end{equation}
	\label{lemma_eq:deform_g}	
\end{lemma}

\begin{proof}
	~~Based on the deformation of the constraint function $g$ along the trajectory $x(t; \zeta,x_0)$ shown in \eqref{eq:deform_g},
	\begin{equation*}
	\frac{d }{dt}g(x(t; \zeta,x_0))= -|\nabla g|(\zeta +\cos \theta(x) ),
	\end{equation*}
we get a proof directly.
\end{proof}

\begin{remark}
Lemma \ref{lemma:deform_g} hints that the optimization progresses faster when an iterate $x(t_k) \in \Theta_{\zeta}$ than vice versa. This result is useful in designing practical algorithms where mechanisms/heuristics can be developed to adjust the stepsize and/or the parameter $\zeta$ so to achieve a faster convergence.
\end{remark}

\subsection{Global convergence to critical points of $f(x)$}
\label{sec:global_conv_critical_points}
\textcolor[rgb]{0,0,0}{
	We show that, as long as $\zeta \in [0,1)$, the present optimization trajectory always converges to a critical point of the objective function.}

\begin{theorem}
Suppose that assumptions (A1),(A3), (A4) hold and $|\nabla g| \neq 0$ in the whole $\mathbb{R}^n$. The trajectory $x(t; \zeta,x_0)$ converges to a connected subset of  critical points of the objective function by $\zeta \in [0,1)$. Especially, if the critical points of the objective function are isolated, then
\begin{equation}\label{convergencetocritical}
\lim_{t\rightarrow T_{\zeta,x_0}^-}x(t; \zeta,x_0)=x_c,~~~~  \forall \zeta \in [0,1),
\end{equation}
where $x_c$ is a critical point of $f$.
\label{theorem:global_conv_critical_point}
\end{theorem}
\begin{proof}
~~First, notice that the present flow may be finite-time convergent\footnote{\textcolor[rgb]{0,0,0}{If $\zeta = 0$, the present system reduces to the normalized gradient flow for unconstrained minimization, which is shown to be finite-time convergent using nontrivial nonsmooth stability analysis \cite[Theorem~8]{cortes2006finite}.}}. The analytical example in section \ref{sec:2d_analytic} shows that the present system can be finite-time convergent: 1) the analytical trajectory has a finite length, and 2) that the speed of the flow has a minimum norm of $1-\zeta$ by the definition of the present system. To overcome this difficulty, we introduce a time-reparameterized Y$-$system:
    \begin{equation}
    \left\{
    \begin{split}
    &\frac{dy}{d\tau}=|\nabla f(y)| \textbf{s}_\zeta(y), \\ 
    & y|_{\tau=0} = x_0.
    \end{split}
    \right.
    \label{Y-system}
    \end{equation}
    The system (\ref{Y-system}) has the same orbit as (\ref{system}) in the subset $E \subset \mathbb{R}^n$. Under our assumptions, the Y-system is uniformly Lipschitz continuous and continuous in $\tau$.  
	\textcolor[rgb]{0,0,0}{	By Picard's existence theorem, it has a unique solution $y(\tau; \zeta,x_0)$ with an infinite existence interval.} 
	
	For $ \zeta \in [0,1) $, consider an integral along $y(\tau; \zeta, x_0)$, \begin{equation*}
	f(y(T; \zeta,x_0))= f(x_0)-\int_0^{T}|\nabla f|^2(1 +\zeta\cos \theta)(y(\tau; \zeta,x_0))d\tau.
	\end{equation*}
	Lemma \ref{lemma_boundedness} ensures
	\begin{equation}
	\int_0^{T}|\nabla f|^2(1 +\zeta\cos \theta)(y(\tau; \zeta,x_0))d\tau= f(x_0)-f(y(T; \zeta,x_0))\leq M 
    \label{eq:bounded_integral}
	\end{equation} 
	with some positive number $M$ independent of $T$. Hence
	\begin{equation}\label{finitint}
	\int_0^{\infty}|\nabla f|^2(y(\tau; \zeta,x_0))d\tau\leq \frac{M}{1-\zeta} \textcolor[rgb]{0,0,0}{< +\infty}.
	\end{equation}
		Notice that
		\begin{equation*}
		\begin{split}
		&\frac{d}{d\tau} | \nabla f (y( \tau; \zeta, x_0)) | = \frac{d}{d \tau} \left( \sum_{j = 1}^{n} \left( \frac{\partial f}{\partial y_j}\right)^2 \right)^{\frac{1}{2}}\\
		&= \frac{1}{2} \left( \sum_{j = 1}^{n} \left( \frac{\partial f}{\partial y_j} \right)^2 \right)^{-\frac{1}{2}}  \sum_{j = 1}^{n} \frac{d}{d \tau} \left( \frac{\partial f}{\partial y_j}\right)^2 \\
		& = \frac{1}{|\nabla f|} \sum_{j = 1}^{n} \frac{\partial f}{\partial y_j} \sum_{k = 1}^{n} \frac{\partial}{\partial y_k} \left( \frac{\partial f}{\partial y_j} \right) \frac{d y_k}{d \tau} \\
		& = \frac{1}{| \nabla f|} \sum_{j = 1, k = 1}^{n} \frac{\partial f}{\partial y_j} \frac{\partial^2 f}{\partial y_j y_k} \frac{d y_k}{d \tau} \\
		& = - \frac{1}{|\nabla f|} \sum_{j = 1, k = 1}^{n} \frac{\partial f}{\partial y_j}\frac{\partial^2 f}{\partial y_j \partial y_k} |\nabla f| \left( \frac{1}{|\nabla f|} \frac{\partial f}{\partial y_k} + \frac{\zeta}{|\nabla g|} \frac{\partial g}{\partial y_k} \right).
		\end{split}
		\end{equation*}
Notice that $\zeta \in [0,1)$, we have
		\begin{equation*}
		\left|\frac{d}{d\tau}|\nabla f(y(\tau; \zeta,x_0))|\right|\leq \sum\limits_{k=1,j=1}^{n} \left|\frac{\partial f}{\partial y_j} \right| \left|\frac{\partial^2 f}{\partial y_k \partial y_j} \right| \left( \frac{1}{|\nabla f|} \left| \frac{\partial f}{\partial y_k} \right| + \frac{1}{|\nabla g|} \left|\frac{\partial g}{\partial y_k}\right| \right).
		\end{equation*}	
		By Cauchy-Schwarz inequality, we have
		\begin{equation*}
		\left|\frac{d}{d\tau}|\nabla f(y(\tau; \zeta,x_0))|\right|\leq |\nabla f|
		\sqrt{\sum\limits_{k=1,j=1}^{n} \left|\frac{\partial^2 f}{\partial y_k \partial y_j} \right|^2} \sqrt{\sum\limits_{k=1}^{n} \left( \frac{1}{|\nabla f|} \left| \frac{\partial f}{\partial y_k} \right| + \frac{1}{|\nabla g|} \left| \frac{\partial g}{\partial y_k} \right| \right)^2} .
		\end{equation*}	
		Further, we have
		\begin{equation*}
		\begin{split}
		\sum_{k=1}^n \left( \frac{1}{|\nabla f|} \left| \frac{\partial f}{\partial y_k} \right| + \frac{1}{|\nabla g|} \left| \frac{\partial g}{\partial y_k} \right| \right)^2 & = \sum_{k=1}^n \left( \frac{1}{|\nabla f|^2} \left| \frac{\partial f}{\partial y_k} \right|^2 + 2 \frac{1}{|\nabla f| |\nabla g|} \left| \frac{\partial f}{\partial y_k} \frac{\partial g}{\partial y_k}\right| +\frac{1}{|\nabla g|^2}\left| \frac{\partial g}{\partial y_k} \right|^2 \right)\\
		& \leq  \sum_{k=1}^n 2 \left( \frac{1}{|\nabla f|^2} \left| \frac{\partial f}{\partial y_k} \right|^2 + \frac{1}{|\nabla g|^2}\left| \frac{\partial g}{\partial y_k} \right|^2   \right) = 4.
		\end{split}
		\end{equation*}
		Therefore, 
		\begin{equation}
		\left|\frac{d}{d\tau}|\nabla f(y(\tau; \zeta,x_0))|\right|\leq 2 
		|\nabla f|
		\sqrt{\sum\limits_{k=1,j=1}^{n} \left|\frac{\partial^2 f}{\partial y_k \partial y_j} \right|^2}.
		\end{equation}
	By assumption (A4) and Lemma \ref{lemma_boundedness}, the right-hand side of above equation is bounded. There is a constant $l$ so that
	\begin{equation*}
	\left|\frac{d}{d\tau}|\nabla f(y(\tau; \zeta,x_0))|\right|\leq l, ~\forall \tau.
	\end{equation*}
	Hence, we have a Lipschitz continuity for $|\nabla f(y(\tau; \zeta,x_0))|$:
	\begin{equation}
	\left||\nabla f(y(\tau^\prime; \zeta,x_0))|-|\nabla f(y(\tau^{\prime\prime}; \zeta,x_0))|\right| \leq l |\tau^\prime-\tau^{\prime\prime}|,~~~~~ \forall \tau^\prime,\tau^{\prime\prime}.
	\label{lipschitz}
	\end{equation}
	Now, we claim
	\begin{equation}
	\lim_{\tau\rightarrow +\infty}|\nabla f(y(\tau;\zeta,x_0))|=0.
	\label{claim}
	\end{equation}
	Otherwise, there is a sequence of $\tau_j\rightarrow +\infty$ and a positive constant $b$ so that
	\begin{equation*}
	|\nabla f(y(\tau_j;\zeta,x_0))|\geq b>0.
	\end{equation*}
	Choosing $\delta=\frac{b}{2l}$, then
	\begin{equation*}
	\begin{split}
	&|\nabla f(y(\tau; \zeta,x_0))|\geq |\nabla f(y(\tau_j; \zeta,x_0))|\\
	&\hskip 3mm -\left||\nabla f(y(\tau_j; \zeta,x_0))|-|\nabla f(y(\tau; \zeta,x_0))|\right| \\
	& \hskip 3mm \geq |\nabla f(y(\tau_j; \zeta,x_0))|-l|\tau_j-\tau|\geq b-\delta l=\frac{b}{2}
	\end{split}
	\end{equation*}
	for any $|\tau_j-\tau|\leq \delta.$ Therefore,
	\begin{equation*}
	\begin{split}
	&\int_0^\infty|\nabla f(y(\tau; \zeta,x_0))|^2d\tau\geq \sum\limits_{j=1}^{\infty} \int_{\tau_j-\delta}^{\tau_j+\delta}|\nabla f(y(\tau; \zeta,x_0))|^2d\tau\\
	&\hskip 3mm \geq \sum\limits_{j=1}^{\infty} \int_{\tau_j-\delta}^{\tau_j+\delta}\frac{b^2}{4}d\tau=\sum\limits_{j=1}^{\infty} \frac{b^3}{4l}=+\infty.
	\end{split}
	\end{equation*}
	This is a contradiction to (\ref{finitint}). Hence
	\begin{equation}
	\lim_{\tau\rightarrow +\infty}|\nabla f(y(\tau;\zeta,x_0))|=0.
	\label{eq:asymptotic_gradient_f}
	\end{equation}
	This means that $y(\tau, \zeta, x_0)$ approaches the connected subset of the critical points. An isolated condition  makes sure that
	\begin{equation}
	\lim_{\tau\rightarrow +\infty}y(\tau; \zeta, x_0)=x_c,
	\label{eq:limit_y}
	\end{equation}
	for some critical point $x_c$ of the objective function. \textcolor[rgb]{0,0,0}{Recall that the system \eqref{Y-system} and \eqref{system} are orbit equivalent in $E$, therefore}
	\begin{equation}
	\lim_{t\rightarrow T_{\zeta,x_0}^-}x(t; \zeta,x_0)=x_c,~~~~~~~~~~ \forall \zeta \in [0,1).
	\end{equation}
	
\end{proof}

\begin{remark}
The global behavior of the system (\ref{system}) looks like a gradient flow. Obviously, as $\zeta = 0$, the system reduces to the normalized gradient flow of the objective function. Specifically, when converging to the same critical point, the present trajectory is homotopic with the gradient descent trajectory for $\zeta \in (0,1)$ by the continuous dependence of the solutions to differential equations on parameters.
\label{remark:critical_point}
\end{remark}

\begin{remark}
The analytical trajectory may find a critical point of constraint $g$ without condition (A3). In practical implementations, we suggest using the gradient descent $-\nabla f(x)$ to escape these critical points.
\end{remark}
\vskip 2mm

\begin{remark}
We note that the assumption (A2) is excluded in Theorem \ref{theorem:global_conv_critical_point}. If there exists some critical point $x_c$ of $f(x)$ in the feasible set $\Omega$, the optimization trajectory may converge to it (also compare to Lemma \ref{lemma:cross_boundary} where the assumption (A2) is included for the analysis). In such a case, obviously, a first-order solution is found for the considered constrained optimization problem. An example is illustrated in figure \ref{fig:2d_nonconvex_obj_quadratic_constraint} for a 2D problem. 
\end{remark}

\subsection{Global convergence to KKT solutions at the constraint boundary} 
In the previous subsection, we have shown that the continuous optimization trajectory converges to a critical point $x_c$ of the objective function $f(x)$, regardless if $x_c$ lies in the feasible set $\Omega$ or not. For cases that $x_c \in \Omega$, this implies that we are able to use the present optimization trajectory to obtain a (first-order) optimal solution. What remains is the question of whether the present method can \textit{find} a (near-optimal) solution at the constraint boundary, i.e., the second case presented in section \ref{sec:nonconvex_example}. In the following, we give an affirmative answer to this question. To give a rigorous presentation, we shall first define what we meant by the wording \textit{``find''}, based on which we establish an error bound of the solution to a KKT point.

\subsubsection{A solution is \textit{found} at the intersection point of the optimization trajectory and the constraint boundary}

As illustrated in figure \ref{fig:2d_linear_constraint_analysis} of the second preliminary example, the optimization trajectories $x(t; \zeta, x_0)$, with different $\zeta$ values, converge to the critical point of the objective function while intersecting the constraint boundary at some time $t$. The larger the $\zeta$ chosen, the closer the intersection point is to the KKT solution. We consider these intersection points as the solutions \textit{found} by the present method to a single constrained problem. A formal definition is given as follows.

\begin{definition} (Solution point $x_\zeta^\sharp$)
Let $\zeta \in [0,1)$ and $x_0 \in \Omega$, and suppose that the optimization trajectory $x(t;\zeta,x_0)$ intersects the constraint boundary at a point $x_\zeta^\sharp = x(t^\sharp; \zeta, x_0)$, where $t_\zeta^\sharp$ is the minimum time for the trajectory to reach the constraint boundary. We say the intersection point $x_\zeta^\sharp$ is a \textit{solution point} obtained by the system \eqref{eq:ds_scop} for the problem \eqref{eq:single_constrained_optimization_problem}. We have
\begin{equation}
\left\{
\begin{split}
&g(x(t^\sharp;\zeta,x_0))=0,\\
& g(x(t;\zeta,x_0))<0,~~ t<t^\sharp.
\end{split}\right.
\label{eq:solution_point}
\end{equation}
\end{definition}

Next, we argue that under the assumption $(A1) - (A4)$, the optimization trajectory $x(t;\zeta, x_0)$, with $\zeta \in [0,1)$ and $x_0 \in \Omega$, always intersects the constraint boundary. Although this statement can readily be made by following Theorem \ref{theorem:global_conv_critical_point}, we give an independent and simple proof for clarity.

\begin{lemma}
	Suppose that assumptions (A1) - (A4) hold, and let $\zeta \in [0,1)$, $x_0 \in \Omega$. Then the trajectory $x(t; \zeta,x_0)$ must go out of the feasible set. Hence, a solution point $x_\zeta^\sharp$ must exist. 
	\label{lemma:cross_boundary}
%
\end{lemma}
\begin{proof}
	~~We prove the lemma by contradiction. Suppose that $x(t; \zeta,x_0)$ stays in the feasible set $\Omega $ for any $t<T_{\zeta, x_0}$, then the vector field $\textbf{s}_{\zeta}$ keeps $C^1$ continuous in a neighborhood of trajectory $x(t; \zeta,x_0)|_{[0,T_{\zeta, x_0})}$, since there is no critical point of $f(x)$ and $g(x)$ in $\Omega$. Lemma \ref{lemma_boundedness} ensures that 
	\begin{equation}
	|\nabla f(x(t; \zeta,x_0))|\geq A, |f(x(t; \zeta,x_0)|\leq B, ~~~ \forall t<T_{\zeta, x_0},
	\label{eq:lowerbound}
	\end{equation}
	with some positive numbers $A$ and $B$. 
	On the other hand, $\mathbf{s}_\zeta(x(t))$ is uniformly Lipschitz continuous in $x(t)$ by Lemma \ref{lemma_lipschitz}. Picard's existence theorem implies that $T_{\zeta, x_0}=+\infty.$
	The integral of (\ref{eq:deform_f}) shows (see also Remark \ref{remark:boundedness})
	$$
	\int_0^{T_{\zeta, x_0}}|\nabla f(x(t; \zeta,x_0))|(1+\zeta \cos \theta)dt = f(x_0)-f(x(T_{\zeta, x_0}; \zeta,x_0)) < \infty.
	$$
    For any $\zeta \in [0,1)$, $1+ \zeta \cos \theta > 0$ holds, therefore
	\begin{equation}
	\int_0^\infty|\nabla f(x(t; \zeta,x_0))|dt<\infty.
	\label{eq:finite_integral}
	\end{equation}
	\eqref{eq:finite_integral} is a contradiction to (\ref{eq:lowerbound}). 	This completes the proof.
\end{proof}

\subsubsection{An upper bound of the solution time $t^\sharp$ relative to $\zeta$} We estimate the solution time $t^\sharp$ that is equivalent for the trajectory $x(t;\zeta,x_0)$ to reach the constraint boundary. An upper bound for $t^\sharp$ is given as follows.
\begin{lemma}
	Suppose that the assumptions (A1)- (A4) hold, and let $\zeta \in [0,1)$, $x_0 \in \Omega$, the time in which the optimization trajectory $x(t;\zeta,x_0)$ reaches the boundary of the feasible set is at most $\frac{C}{1-\zeta}$ with $C$ independent of $\zeta$. 
	\label{lemma:time_to_boundary}
\end{lemma}
\begin{proof}
	~~By Lemma \ref{lemma:cross_boundary}, the trajectory $x(t;\zeta,x_0)$ must go out of the feasible set. Let $T_{\zeta, x_0}^\sharp$ be the time in which the trajectory first reaches the boundary of the feasible set, then (\ref{eq:lowerbound}) holds for $0<t<T_{\zeta, x_0}^\sharp$. By \eqref{eq:deform_f},
	\begin{equation}
	\int_0^{T_{\zeta, x_0}^\sharp}|\nabla f|(1+\zeta \cos \theta)dt=f(x_0)-f(x(T_{\zeta, x_0}^\sharp; \zeta,x_0))\leq 2B.
	\label{eq:df_bound}
	\end{equation}	
	Hence
	\begin{equation}
	T_{\zeta, x_0}^\sharp A(1-\zeta)\leq 2B,
	\end{equation}
	which implies
	\begin{equation}
	T_{\zeta, x_0}^\sharp\leq \frac{2B}{A(1-\zeta)},
	\end{equation}
	so that $T_{\zeta, x_0}^\sharp \leq \frac{C}{1-\zeta}$ holds. This ends the proof.
\end{proof}

Lemma \ref{lemma:time_to_boundary} implies that the system \eqref{eq:ds_scop} finds a solution $x_\zeta^\sharp$ in $t^\sharp \leq O(\frac{1}{1-\zeta})$.

\subsubsection{An upper bound of the solution error relative to $\zeta$}
Lastly, we shall give an error measure of a solution (intersection) point $x_\zeta^\sharp$ to a KKT solution $x^\star$. The KKT conditions for the problem \eqref{eq:single_constrained_optimization_problem} write
\begin{subequations}
\begin{align}
\nabla f(x^\star) + \lambda^\star \nabla g(x^\star) &= 0, \label{eq:kkt_x_sharp1}\\
\lambda^\star = \frac{|\nabla f(x^\star)|}{|\nabla g(x^\star)|} &\geq 0, \label{eq:kkt_x_sharp2}\\
g(x^\star) &= 0.
\label{eq:kkt_x_sharp3}
\end{align}
\end{subequations}
It is obvious that the conditions \eqref{eq:kkt_x_sharp1} and \eqref{eq:kkt_x_sharp2} are equivalent to the normalized centrality condition \eqref{eq:normalized_central_path_condition_single_constraint}, repeated as the follows for convenience,
\begin{equation*}
\frac{\nabla f(x)}{|\nabla f(x)|} + \frac{\nabla g(x)}{|\nabla g(x)|}  = 0,
\end{equation*}
In addition, \eqref{eq:kkt_x_sharp3} is satisfied with \eqref{eq:solution_point}. Following the normalized central path condition \eqref{eq:normalized_central_path_condition_single_constraint}, a proper residual for an intersection point $x_\zeta^\sharp$ can be defined using the L2-norm,
\begin{equation}
r_{L2}(x_\zeta^\sharp) = \left| \frac{\nabla f(x_\zeta^\sharp)}{| \nabla f(x_\zeta^\sharp) |} + \frac{\nabla g(x_\zeta^\sharp)}{| \nabla g(x_\zeta^\sharp)|} \right|.
\label{eq:dual_residual}
\end{equation}
Based on \eqref{eq:dual_residual}, we give a formal definition for an error measure $\epsilon$.
\vskip 2mm
\begin{definition} [Solution error]
	An error measure $\epsilon$ for a solution point $x_\zeta^\sharp$ is defined as the L2-norm of the residual \eqref{eq:dual_residual},
	\begin{equation}
	\epsilon(x^\sharp) = r_{L2}(x^\sharp).
	\label{eq:epsilon_dr}
	\end{equation}
	\label{def: epsilon}
\end{definition}
In the following, we show that the error $\epsilon$ of the method is upper bounded relative to the parameter $\zeta$.

\begin{theorem}
Suppose that assumptions (A1)-(A4) hold, and let $\zeta \in [0,1)$, $x_0 \in \Omega$. Then, the solution error of $x_\zeta^\sharp$ is upper bounded in relative to the parameter $\zeta$,
  \begin{equation}
      \epsilon(x_\zeta^\sharp) \leq  \sqrt{2(1-\zeta)}.
\label{eq:solution_error}
  \end{equation}
\label{theorem:global_error_bound}
\end{theorem}
\begin{proof}
~~By Lemma \ref{lemma:cross_boundary}, the optimization trajectory $x(t;\zeta, x_0)$ must intersect the constraint boundary when initialized in the feasible set $\Omega$. By \eqref{eq:solution_point}, we have at the intersection point $x(t^\sharp; \zeta, x_0)$,
\begin{equation}
\frac{d }{dt}g(x_\zeta^\sharp)\geq 0.
\end{equation}
By (\ref{eq:deform_g}),
\begin{equation}
 |\nabla g(x_\zeta^\sharp)|(\zeta +\cos \theta(x_\zeta^\sharp))\leq 0.
\end{equation}
Thus,
\begin{equation}
 \cos \theta(x_\zeta^\sharp)\leq -\zeta.
 \label{eq:cos_theta_zeta}
\end{equation}
By Definition \ref{def: epsilon}, we have
\begin{equation}
\epsilon^2(x_\zeta^\sharp) = 2(1+\cos \theta(x_\zeta^\sharp)).
\end{equation}
Therefore,
\begin{equation*}
\epsilon(x_\zeta^\sharp) \leq \sqrt{2(1-\zeta)}.
\end{equation*}
Our proof is thus complete.
\end{proof}

\begin{remark}
By Lemma \ref{lemma:time_to_boundary} and Theorem \ref{theorem:global_error_bound}, we can conclude that the time complexity for the optimization trajectory to achieve an approximate solution is
\begin{equation}
| x_\zeta^\sharp - x^\star| \leq O(\frac{1}{\sqrt{t}}).
\end{equation}
The result matches the ergodic convergence rate of first-order methods for general optimization problems.
\end{remark}

\section{The method for multiple constraints: a formulation based on the logarithmic barrier function}
\label{sec:multiple_constraint}
Recall the logarithmic barrier function for the problem (\ref{eq:Optimization_Problem}),
$$\Phi(x) = - \sum_{i = 1}^{m} \log (-g_i(x)), ~ i = 1,...,m.$$
We first notice that the barrier function grows without bound if $g_i(x) \rightarrow 0^-$ and it is not differentiable at the boundary of the feasible set. To overcome this difficulty, we consider a subset of the feasible set,
\begin{equation}
\Omega_M=\{x: \Phi(x) \leq M \},
\end{equation}
where M is a sufficiently large positive number so that $\Omega_M$ approximates the original feasible set $\Omega$.
The barrier function is twice continuously differentiable in $\Omega_M$ and thus fulfills the assumption (A4). 

Set $$ G_M(x) = \Phi(x) - M, ~ x \in \Omega_M. $$
The original problem (\ref{eq:Optimization_Problem}) can be approximately reformulated as
\begin{equation}
\begin{split}
&\textnormal{minimize}~~~  f(x), \\
&\textnormal{subject to}~~ G_M(x) \leq 0, \\
\end{split}
\label{eq:Log_barrier_approximation}
\end{equation}
so that the present system (\ref{eq:ds_scop}) may be applied. 

Notice that the global behavior shown in section \ref{sec: global_behavior} is based on the assumption that $\nabla g \neq 0$ in the feasible set $\Omega$. In general, this may not always be the case for the function $G_M(x), \forall x \in \Omega_M$. To escape the points where $\nabla \Phi(x) = 0$, we suggest using the gradient descent $-\nabla f(x)$. Thus we get a dynamical system for solving the approximate problem \eqref{eq:Log_barrier_approximation},
\begin{equation}
\tag{DSM}
\frac{d x}{d t}=\left\{
\begin{split}
&-\frac{\nabla f(x)}{|\nabla f(x)|}
-\zeta\frac{\nabla \Phi(x)}{|\nabla \Phi(x)|},  ~~~~~&\text{if}~~\nabla \Phi(x) \neq 0;\\
& -\nabla f(x)    &\text{if}~~\nabla \Phi(x) = 0.
\end{split}\right.
\label{eq:ds_mcop}
\end{equation}
Note that in computational practice, $\nabla \Phi = 0$ rarely occurs. The second ODE of the above system serves as a safeguard for the present method.
\vskip 2mm

\begin{corollary}
Suppose that the assumptions (A1) (A2) and (A4) hold, and $\nabla g_i(x) \neq 0, i = 1,...,m$ in the feasible subset $\Omega_M$. Let $x_\zeta^\sharp$ be the solution point of the system \eqref{eq:ds_mcop} for the approximate problem \eqref{eq:Log_barrier_approximation}, then
\begin{equation}
x_\zeta^\sharp\in \{x:\Phi(x) = M , \cos\theta(x) \leq -\zeta \}.
\end{equation}
And the solution error to a KKT solution of the problem \eqref{eq:Log_barrier_approximation} is
$$\epsilon \leq O(\sqrt{1-\zeta}).$$
\end{corollary}

\begin{proof}
~~A similar method as used in Theorem \ref{theorem:global_error_bound} gives the proof.	
\end{proof}

To conclude, using the barrier function formulation for problem (\ref{eq:Optimization_Problem}), the resulting trajectory achieves an approximated local solution that locates on the boundary of the subset $\Omega_M$. Notice, too, that the system \eqref{eq:ds_mcop} does not depend on the choice of $M$, and $\Omega_M$ exhaust $\Omega$, i.e., $\Omega = \cup_{M>0} \Omega_M$. This means that the resulting trajectory keeps approaching the boundary of the original feasible set $\Omega$ by crossing the boundary of any subset $\Omega_M$ dependent on $M$. Therefore, it eases the practical implementation of the method since no extra parameter $M$ needs to be defined for the stopping criterion, but rather, a check for each constraint violation would be sufficient.

\section{Algorithms and preliminary computational experiments}
\label{sec:algorithms}
We present two computational implementations of GDAM and show their preliminary computational tests on common benchmarks for both convex and nonconvex constrained optimizations. 

\subsection{Computational algorithms} \label{sec:implementation} 
When implementing the present method for the problem \eqref{eq:Optimization_Problem}, the canonical first-order optimization procedure is given as below,
\begin{equation}
x_{k+1} = x_k + \alpha_k \mathbf{s}_\zeta(x_k),
\label{eq:GD_trajectory}
\end{equation}
with 
\begin{equation*}
\mathbf{s}_\zeta(x_k)=\left\{
\begin{split}
&-\frac{\nabla f(x_k)}{|\nabla f(x_k)|}
-\zeta\frac{\nabla \Phi(x_k)}{|\nabla \Phi(x_k)|},  ~~~~~&\text{if}~~\nabla \Phi(x_k) \neq 0;\\
& -\nabla f(x_k)    ~~~~~~~~~~~~~~~~~~~~&\text{if}~~\nabla \Phi(x_k) = 0.
\end{split}\right.
\end{equation*}

In this work, we design computational algorithms for problems whose minimizers are non-trivially located at the constraint boundary, as is the case for our targeted applications, shape optimization and sensor network localization. 

\subsubsection{Fix-length stepsize}
Obviously, the choice of $\alpha_k$ is essential to the computational performance of the method, and, certainly, there are many ways that can satisfy the purpose of the convergence. In this work, we propose a fixed-length step size rule for the method, i.e.,
\begin{equation}
    \alpha_k = \frac{\beta}{| s_\zeta (x_k)|},
\label{eq:fixed_length}
\end{equation}
where $\beta$ is some positive parameter. As a result,
\begin{equation}
| x_{k-1} - x_k | = \beta.
\end{equation}
The step size rule can be considered aggressive, since at a central point $x_c$, as $\zeta \rightarrow 1^-$, $| s_\zeta(x_c) | \rightarrow 0$. A fixed-length step size \eqref{eq:fixed_length} aggressively re-scales the step size using the reciprocal of $| s_\zeta(x) |$, resulting in a finite variable update. A fix-length stepsize is also used in subgradient methods to deal with nonsmoothness at the solution by dynamically adjusting $\beta$ towards zero, see, e.g., \cite[Proposition 3.2.7]{bertsekas2015convex}. For constrained problems, when more than one inequality constraints are active at a solution, the solution vertex itself can be nonsmooth and contributes to the ill-conditioning of the logarithmic barrier function at its close neighborhood. A fix-length stepsize bounds the range of an update and can thus be considered a trust-region-like approach that stabilizes the iterative process. Finally, various line search methods can be added to dynamically adjust the parameter $\beta$. 

\subsubsection{Vanilla implementation}
A vanilla implementation based on the fixed-length stepsize approach is summarized in the pseudocode in Algorithm \ref{alg:1}. For constant parameter $\beta$, we terminate the optimization whenever a constraint is violated. 

\begin{algorithm}
\caption{Vanilla GDAM}\label{alg:1}
\begin{algorithmic}[1]
\State \textbf{Initialization:} $x_0$, $\zeta$, $\beta$, $\beta_{min}$, $k = 0$
\WHILE{Stopping criteria is not met}
\IF{$\nabla \Phi(x_k) \neq 0$}
\Statex $ ~~~~\mathbf{s}_\zeta(x_k) \leftarrow - \frac{\nabla f(x_k)}{| \nabla f(x_k) |} - \zeta \frac{\nabla \Phi(x_k)}{ | \nabla \Phi(x_k) |} $ 
\ELSE 
\Statex		$ ~~~~\mathbf{s}_\zeta(x_k) \leftarrow - \nabla f(x_k)$
\ENDIF
\State \textbf{Update:} 
$\alpha_k \leftarrow \frac{\alpha}{| s_\zeta (y_k)|}$,
$x_{k+1}  \leftarrow x_k + \alpha_k  s_\zeta (x_k),$ 
$k \leftarrow k+1$
\State \textbf{(Line search to update $\beta$)}
\ENDWHILE
\end{algorithmic}
\end{algorithm}
When an additional line search is equipped, the stopping criteria is set as to simultaneously fulfill the following two conditions:
\begin{itemize}
    \item [1.] Any constraint is violated;
    \item [2.] $\beta < \beta_{min}$.
\end{itemize}
For example, the algorithm can be equipped with a simple backtracking line search to improve the performance, i.e., we reduce the stepsize $\beta$ by a scaling factor $\tau$ whenever a constraint is violated until $\beta < \beta_{min}$.

\subsubsection{Accelerated implementation}
A prominent approach to improve the convergence rate of the gradient descent method is Nesterov's Accelerated Gradient (NAG) method \cite{nesterov1983method}. We apply NAG to the present method by replacing the gradient descent update step with a GDAM update step. We summarize the pseudocode in Algorithm \ref{alg:2}. Notice that compared to the vanilla implementation, a momentum parameter $m$ and a sequence of auxiliary variables $\{y_k\}$ are added. 

\begin{algorithm}
\caption{Accelerated GDAM}\label{alg:2}
\begin{algorithmic}[1]
\State \textbf{Initialization:} $y_0$, $\beta$, $\beta_{min}$, $m$, $k = 0$ 
\State \textbf{Compute:} 
$\alpha_0 = \frac{\alpha}{| s_\zeta (y_0)|}$,
 $x_0 \leftarrow y_0 + \alpha_0  s_\zeta (y_0)$, 
$y_1 \leftarrow x_0$, $k = 1$
\WHILE{Stopping criteria is not met}
\IF{$\nabla \Phi(x_k) \neq 0$}
\Statex $ ~~~~\mathbf{s}_\zeta(y_k) \leftarrow - \frac{\nabla f(y_k)}{| \nabla f(y_k) |} - \zeta \frac{\nabla \Phi(y_k)}{ | \nabla \Phi(y_k) |} $ 
\ELSE 
\Statex		$ ~~~~\mathbf{s}_\zeta(y_k) \leftarrow - \nabla f(y_k)$
\ENDIF
\State \textbf{Update:} 
$\alpha_k \leftarrow \frac{\beta}{| s_\zeta (y_k)|}$,
$x_k  \leftarrow y_k + \alpha_k  s_\zeta (y_k),$
$y_{k+1} \leftarrow x_k + m (x_k - x_{k-1}),$
$k \leftarrow k+1$
\State \textbf{(Line search to update $\beta$)}
\ENDWHILE
\end{algorithmic}
\end{algorithm}
\bigskip

The following shows preliminary computational experiments of the proposed vanilla and accelerated GDAM algorithms on common benchmark tests. Computational results show that the present implementations are robust in finding near-optimal solutions, which can be desirable in practical engineering applications.

\subsection{Experiments on IEEE CEC 2006 tests with the vanilla implementation}
We first show numerical experiments for the inequality constrained problems presented at the EA competition at the 2006 IEEE Congress on Evolutionary Computation \cite{liang2006problem}. These benchmark tests are widely used in the community of evolutionary algorithms. We choose them as test examples for three reasons. First, they are well defined constrained optimization problems and have different characteristics \cite{rao2016jaya}. Second, they are nontrivial to solve with first-order methods. Third, the relatively simple formulation of the optimization problems allows us to gain deeper insight into the numerical behavior of the present method. We choose the inequality constrained optimization problems for the experimentation. The problem $G12$ is excluded because it has a feasible set consisting of $9^3$ disjointed spheres. For the problem $G24$, which has a feasible set consisting of two disconnected sub-regions, we choose initial designs in the sub-region that contains the reported optimal solution.

We conduct numerical experiments for the vanilla implementation (Algorithm \ref{alg:1}) with $\zeta = 0.98$ and tuned fixed-length stepsize parameter $\beta$. All bound constraints are treated as inequalities. Random initializations that are away from the reported optimal solution are selected in the feasible sets. We evaluate the error in the objective function value by
\begin{equation}
    \text{err} = \frac{ | f(\mathbf{x}^\star) - f(\mathbf{x_\zeta}) |}{1 + | f(\mathbf{x}^\star) | },
\end{equation}
where $\mathbf{x}^\star$ is the reported global optimum. As shown in table \ref{tb:stepsize_0.98}, apart from the problem G19, vanilla GDAM finds solutions with objective function value errors less than 2e-2 compared to the reported global optima. This is somewhat surprising because GDAM is a local search algorithm. More accurate results can be obtained when we choose shorter step sizes and larger parameters $\zeta$. For the problem G19, the reported optimal solution $ \mathbf{x}^\star = $ (1.6699e-17, 3.9637e-16, 3.9459, 1.060e-16, 3.2831, 9.9999, 1.1283e-17, 1.2026e-17, 2.5071e-15, 2.2462e-15, 0.3708, 0.2785, 0.5238, 0.3886, 0.2982) is not achievable with the present fixed-length stepsize rule with a reasonable constant parameter $\beta$. For the problems G01, G04, G06, G07, G08, G24, we find sub-optimal solutions that are close to the reported optimal designs. For the problems G09, G10, G18, G19, sub-optimal solutions close to local minima are found. 

\begin{table} [h]
	\begin{center}
	\caption{Results of Algorithm \ref{alg:1} with parameter $\zeta = 0.98$ and tuned fixed step sizes}
	\captionsetup{justification=centering}
	\begin{tabular}{|c | c | c | c | c | c |}
		\hline
		Prob. &  step size&  iters. & f($\mathbf{x}^\star$) & f($\mathbf{x_\zeta}$) & obj. error   \\ [0.5ex]
		\hline\hline
		G01  &  0.002 & 2362 & -15 & -14.7215 & 1.74e-2\\  
		\hline
		G04  &  0.2   & 136  & -3.0665e+4 & -3.0657e+4  & 2.61e-4  \\
		\hline
		G06 &  0.002 & 4826 & -6.9618e+3 & -6.8371e+3 & 1.79e-2\\
		\hline
		G07 &  0.0027 & 3009 & 24.3062 & 24.7876 & 1.90e-2 \\
		\hline
		G08 &  0.01 & 66 & -9.5825e-2& -9.5063e-2 & 6.95e-4\\
		\hline
		G09 &  0.05 & 120 & 6.8063e+2 & 6.9238e+2 & 1.72e-2\\
		\hline		
		G10 &  0.35 & 5319 & 7.0492e+3 & 7.1898e+3 & 1.99e-2 \\
		\hline	
		G18 &  0.01 & 257 & -0.8660 & -0.8546 & 6.10e-3 \\
		\hline	
		G19 &  0.05 & 294 & 32.6556 & 2.7120e+2 & 7.09 \\
		\hline	
		G24 &  0.02 & 268 & -5.5080 &  -5.4147 & 1.43e-2 \\
		\hline			
	\end{tabular}
	\label{tb:stepsize_0.98}	
	\end{center}
\end{table}

\subsection{Experiments on Maros and Meszaros QPs with the accelerated implementation}
The Maros and Meszaros test set \cite{MarosMeszaros} is a collection of hard convex quadratic programming problems from a variety of sources. The problems included are of the following form
\begin{equation}
\begin{split}
&\min_x ~~~  \frac{1}{2} x^T Hx + c^T x + c_0, \\
&\text{~s.t.} ~~~~  Ax = b, \\
& ~~~~~~~~~  l \leq x \leq u,
\label{eq:maros_meszaros_qp}
\end{split}
\end{equation}
where $x \in \mathbb{R}^n$, $H \in \mathbb{R}^{n \times n}$ is a symmetric positive definite matrix, $A \in \mathbb{R}^{m \times n}$ is a sparse matrix, $c, c_0 \in \mathbb{R}^n$, $ b \in \mathbb{R}^m$, and $l, u \in \mathbb{R}^n$. Note that components of $l$ and $u$ may be infinite, and $ | c_0  |  < +\infty $.

The algorithms \ref{alg:1} and \ref{alg:2} can easily be extended to handle linear equality constraints by the gradient projection method. Suppose that $A$ is nondegenerate, following Rosen's method \cite{rosen1960gradient}, the projection matrix for the subspace spanned by the linear equalities $Ax = b$ writes
\begin{equation}
    P_A = I_n - A^T ( A A^T )^{-1} A.
\end{equation}
The projected vector $v_p$ of any vector $v \in \mathbb{R}^n$ reads
\begin{equation}
    v_p = P_A v.
\end{equation}
Then, the projected GDAM search direction writes,
\begin{equation}
    \mathbf{s}_{\zeta, P}(x_k) = - \frac{P_A \nabla f(x_k)}{| P_A \nabla f(x_k) |} - \zeta \frac{P_A \nabla \Phi(x_k)}{ | P_A \nabla \Phi(x_k) |}.
\end{equation}

Practically, we can precompute and store the projection matrix $P_A$. For large-scale and sparse problems ($n$ and $m$ are large, and $A$ is sparse), we may make use of Cholesky factorization for $AA^T$ for more efficient computation. More advanced methods for treating linear equalities, such as \cite{brust2022large}, may be applied to the present method. 

We use a subset of the Maros and Meszaros instances for the experiments. These instances were carefully chosen to cover a variety of problem models for performance comparisons among dedicated QP solvers \cite{mihic2021managing}. We evaluate the error in the objective function value by
\begin{equation}
    \text{err} = \frac{ | \text{obj}_{\text{best}} - \text{obj} |}{1 + |\text{obj}_{\text{best}} | },
\end{equation}
where $\text{obj}_{\text{best}}$ is the smallest objective function value of the solvers. In table \ref{tb:maros_meszaros}, we show the results of the present method and compare it with well-established solvers Matlab quadprog, Gurobi \cite{gurobi100}, and OSQP \cite{stellato2020osqp}. Both Matlab quadprog and Gurobi implement state-of-the-art second-order solvers for convex QPs, and OSQP is a first-order convex QP solver based on the alternating direction method of multipliers \cite{boyd2011distributed}. For Matlab quadprog and Gurobi, we use the default settings. For OSQP, we set the absolute and relative tolerance to $1e-4$ for more accurate solutions (default settings are both $1e-3$). For GDAM, $\zeta = 0.999$ and a backtracking line search with reduction parameter $\tau = 0.3$ are chosen. We set the maximum number of iterations to 10000 for both OSQP and GDAM. All experiments are run on a Linux workstation with a 2.70 GHz 6-core AMD Ryzen 5 5600U processor and 16 GB RAM.

\begin{table}[h]
\caption{Large Maros and Meszaros instances (Matlab format credit to \cite{mihic2021racqp})}
  \resizebox{0.7\textwidth}{!}{\begin{minipage}{\textwidth}
\begin{tabular}{|c|c|c|c|c|c|c|c|c|c|c|c|c|c|c|}
\hline
 Problem &  \multicolumn{4}{|c|}{Runtime [s]} & \multicolumn{4}{|c}{Num. iters.} & \multicolumn{4}{|c|}{Obj. error}\\
{}   & Matlab   & Gurobi    & OSQP   & GDAM & Matlab   & Gurobi    & OSQP   & GDAM & Matlab   & Gurobi    & OSQP   & GDAM\\
\hline
AUG2D      &  0.37 & 0.02   & 0.03  & 0.27 & 1   & 2 & 25 & 300 & 7.06e-16   & opt\footnote{A feasible solution that has the lowest objective function value is chosen as reference and is considered ``optimal''.} & 5.84e-6 & 2.33e-7\\
AUG2DC     &  0.11 & 0.02   & 0.03  & 0.33 & 1   & 2 & 25 & 350 & opt   & 2.62e-16 & 5.81e-6 & 1.29e-7 \\
AUG2DCQP   &  0.95 & 0.06   & 0.59  & 1.46 & 10   & 16 & 1075 & 1050 & 4.89e-12   & opt & 7.30e-5 & 3.36e-2\\
AUG2DQP    &  0.90 & 0.06   & 0.55  & 1.53 & 10   & 16 & 975 & 1100 & 7.17e-12   & opt & 1.98e-4 & 3.34e-2 \\
AUG3D      &  0.09 & 0.004   & 0.008  & 0.05 & 1   & 2 & 25 & 400 & 5.32e-15   & opt & 2.02e-10 & 1.03e-4 \\
AUG3DC     &  0.01 & 0.006   & 0.009  & 0.05 & 1   & 2 & 50 & 350 & opt   & 1.39e-15 & 6.69e-11 & 4.24e-5 \\
AUG3DCQP   &  0.11 & 0.02   & 0.009  & 0.16 & 8   & 15 & 50 & 640 & 1.21e-8   & opt & 3.09e-6 & 8.61e-5\\
AUG3DQP    &  0.12 & 0.02   & 0.01  & 0.13 & 7   & 15 & 75 & 500 & 2.93e-7   & opt & 1.32e-9 & 4.67e-4 \\
BOYD1     &  24.56 & 0.36   & 9.76  & 28.70 & 26   & 26 & 2750 & 10000 & 4.31e-12   & opt & 5.26 & 1.97e-2\\
CONT-050  &  0.29 & 0.02   & 1.01  & 0.32 & 11   & 15 & 5150 & 900 & 2.83e-6   & opt & 1.24e-4 & 2.18e-4\\
CONT-100  &  0.50 & 0.11   & 12.45  & 2.49 & 4   & 15 & 10000 & 1450 & 2.48e-2   & opt & 2.80e-4 & 1.47e-3\\
CONT-101  &  0.31 & 0.09   & 14.67  & 0.83 & 3   & 11 & 10000 & 452 & 0.69   & opt & 1.44e-2 & 5.97e-7\\
CONT-200  &  3.76 & 0.75   & 119.59 & 9.89 & 4   & 20 & 10000 & 750 & 0.30   & opt & 0.46 & 3.91e-3\\
CONT-201  &  1.93 & 0.52   & 112.91  & 5.21 & 3   & 14 & 10000 & 397 & 0.67   & opt & 0.69 & 1.06e-6\\
CONT-300  &  5.13 & 1.32   & 295.89  & 22.02 & 3   & 14 & 10000 & 715 & 0.67   & opt & 0.52 & 7.56e-6\\
CVXQP1\_L &  11.48 & 12.66   & 61.25  & 0.56 & 13   & 59 & 9075 & 800 & opt   & 4.63e-8 & 2.23e-5 & 3.24e-5 \\
CVXQP2\_L &  4.81 & 1.23   & 5.84  & 0.27 & 13   & 13 & 425 & 850 & 7.66e-12 & opt & 3.63e-9 & 3.12e-5 \\
CVXQP3\_L &  17.06 & 8.99   & 29.56  & 0.83 & 15   & 51 & 3225 & 800 & opt   & 8.52e-6 & 4.62e-5 & 2.76e-5\\
DTOC3 &  0.03 & 0.02   & 0.03  & 0.19 & 1   & 2 & 50 & 400 & opt   & 1.37e-13 & 0.18 & 8.73e-5\\
HUES-MOD  &  0.25 & 0.01   & 0.03  & 0.17 & 11   & 14 & 175 & 650 & 1.46e-11   & opt & $\infty$ & 2.03e-5\\
HUESTIS   &  0.26 & 0.01   & 0.01  & 0.16 & 11   & 15 & 50 & 600 & opt   & 5.78e-15 & $\infty$ & 1.76e-5\\
STCQP1    &  4.06 & 0.09   & 0.04  & 4.17 & 9   & 11 & 150 & 650 & 4.42e-10   & opt & 2.78e-8 & 9.99e-5\\
STCQP2    &  0.59 & 0.16   & 0.06  & 0.11 & 9   & 11 & 125 & 500 & 4.55e-10   & opt & 9.99e-7 & 1.72e-4\\
UBH1    &  0.16 & 0.01   & 0.03  & 1.78 & 3   & 7 & 50 & 2300 & 6.72e-11  & opt & 1.21 & 8.68\\
\hline
\end{tabular}
\end{minipage}
}
\label{tb:maros_meszaros}
\end{table}

As shown in table \ref{tb:maros_meszaros}, accelerated GDAM finds near-optimal solutions for all instances except UBH1, which demonstrates its robustness. Gurobi is the fastest and the most robust solver in general. GDAM is more efficient in solving the instances CVXQP*. By setting smaller tolerances and a larger number of maximum iterations, OSQP finds more accurate solutions for the instances BOYD1, CONT*, and DTOC3. We note that no presolve and preconditioning strategies have been implemented for GDAM, unlike other mature QP solvers, since the goal here is to show the robustness of the present algorithm. On the other hand, it leaves room for improvements in computational efficiency.

\section{Application to Shape Optimization}
\label{sec:shapopt}
We show applications of vanilla GDAM (Algorithm \ref{alg:1}) to shape optimization problems arising in Computational Mechanics. We begin by briefly introducing the shape optimization framework used in this work and then show the optimization results using GDAM.

\subsection{Vertex Morphing}

The shape optimization problems we aim to solve are almost always ill-posed due to the very large discrete design space based on the finite element mesh. To tackle the problem, we use the Vertex Morphing method (VM) introduced in \cite{Hojjat}.

The idea of the Vertex Morphing is to control the discrete surface coordinates $\textbf{x} = [x_1, x_2,..., x_n]^T$ with design controls $\textbf{p} = [p_1, p_2,..., p_n]^T$, filtered by a filter function. The explicit filtering used in VM is the convolution of the coordinate field $x$ with a kernel. In the discretized system, it is a matrix-vector multiplication that is the summation of nodal contributions that are weighted with the kernel function. VM distinguishes itself with standard explicit filtering, in which it applies the filtering process twice.

First, the so-called forward mapping step that uses the linear filtering matrix $\textbf{R}$ is defined as follows:

\begin{equation}
x_i = R_{ij} p_j.
\end{equation}

Similarly, the change of the control $\delta \textbf{p}$ is mapped onto the change of the design configuration $\delta \textbf{x}$

\begin{equation} \label{eq:forward}
\delta x_i = R_{ij} \delta p_{j}.
\end{equation}

The design controls are the control variables of the gradient-based optimization. The number of variables is equivalent to the number of surface coordinates. Following the chain rule of differentiation, the sensitivities of a response function $\Psi$ with respect to the discretized geometry $\textbf{x}$ are backward mapped to the design control using the adjoint or backward mapping matrix $\textbf{R}^*$, with $\textbf{R}^*=\textbf{R}^T$ for regular grids,

\begin{equation} \label{eq:backward}
\frac{d \Psi}{dp_i} = \frac{d \Psi}{d x_j} \frac{d x_j}{d p_i} = R_{ji} \frac{d \Psi}{d x_j}.
\end{equation}

Equations (\ref{eq:forward}) and (\ref{eq:backward}) can be used in a gradient descent framework, ensuring smooth shape updates in each iteration.

\subsection{Academic example}
First, we show an academic convex problem. The objective is to maximize the volume of a small sphere. This small sphere is located inside a bigger sphere that acts as the geometric constraint. The shape of both spheres are represented with finite element meshes \cite{zienkiewicz1977finite}. The optimization problem writes
\begin{equation}
\begin{split}
& \textnormal{minimize} ~~~   -V(x), \\
& \textnormal{subject to} ~~~  g_i(x) \leq 0, ~ i = 1,...,m,\\
\end{split}
\end{equation}
where $V(x)$ is the volume function, $g_i(x)$ is a point-wise defined geometric constraint for the $i$-th design node, $m$ is the number of nodes of the design mesh, and $x \in \mathbb{R}^{3m}$ is the field of nodal coordinates of the design sphere mesh. The number of nodes of the small sphere (design sphere) is 19897. Thus, the total number of design variables is 59691 and the total number of constraints is 19897. We use the logarithmic barrier function for multiple constraints and choose the parameter $\zeta = 0.95$. In figure \ref{fig:sphere_iters}, the shape variation process with depicted iterations is shown. Initially, the design sphere is located close to the boundary of the constraint sphere. During the shape variation process, it moves towards the center of the constraint sphere while adapting its shape at each iteration. One could recognize easily that the central path of this optimization problem is being approached and followed until the solution is found when constraints become active.

\begin{figure}[h]
\centering
\subfloat[Iteration 1]{%
\resizebox*{3.7cm}{!}{\includegraphics{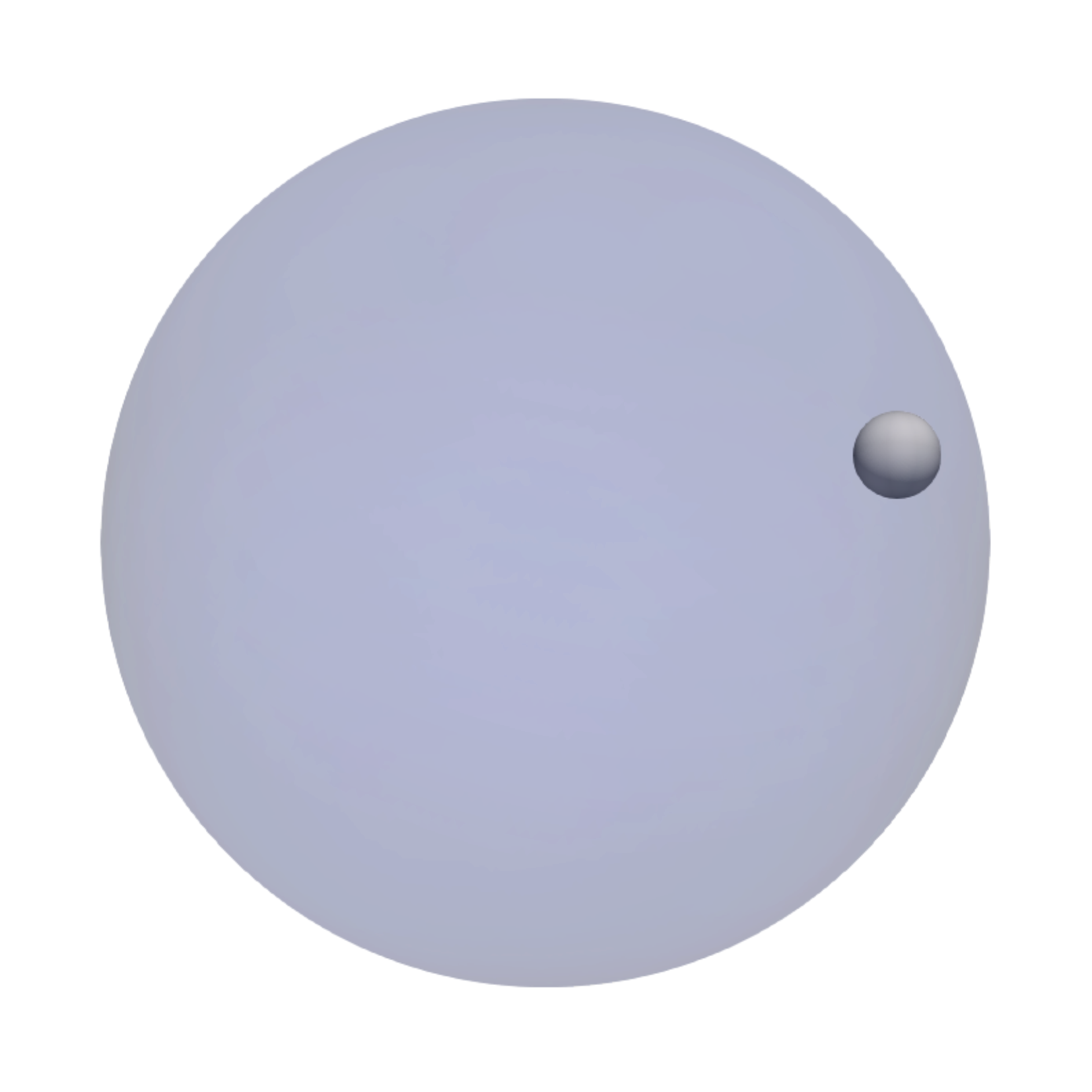}}}
\subfloat[Iteration 15]{%
\resizebox*{3.7cm}{!}{\includegraphics{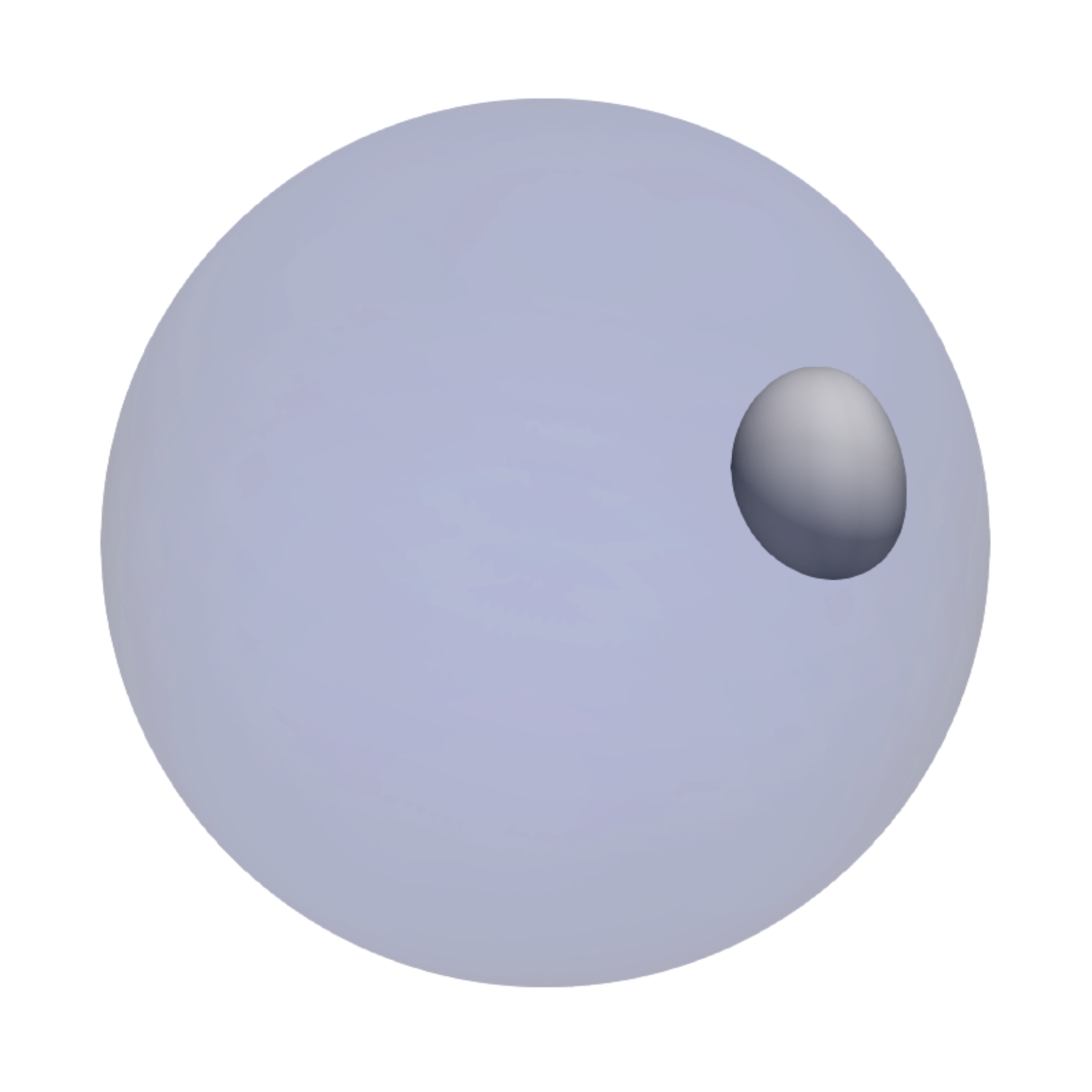}}}
\subfloat[Iteration 30]{%
\resizebox*{3.7cm}{!}{\includegraphics{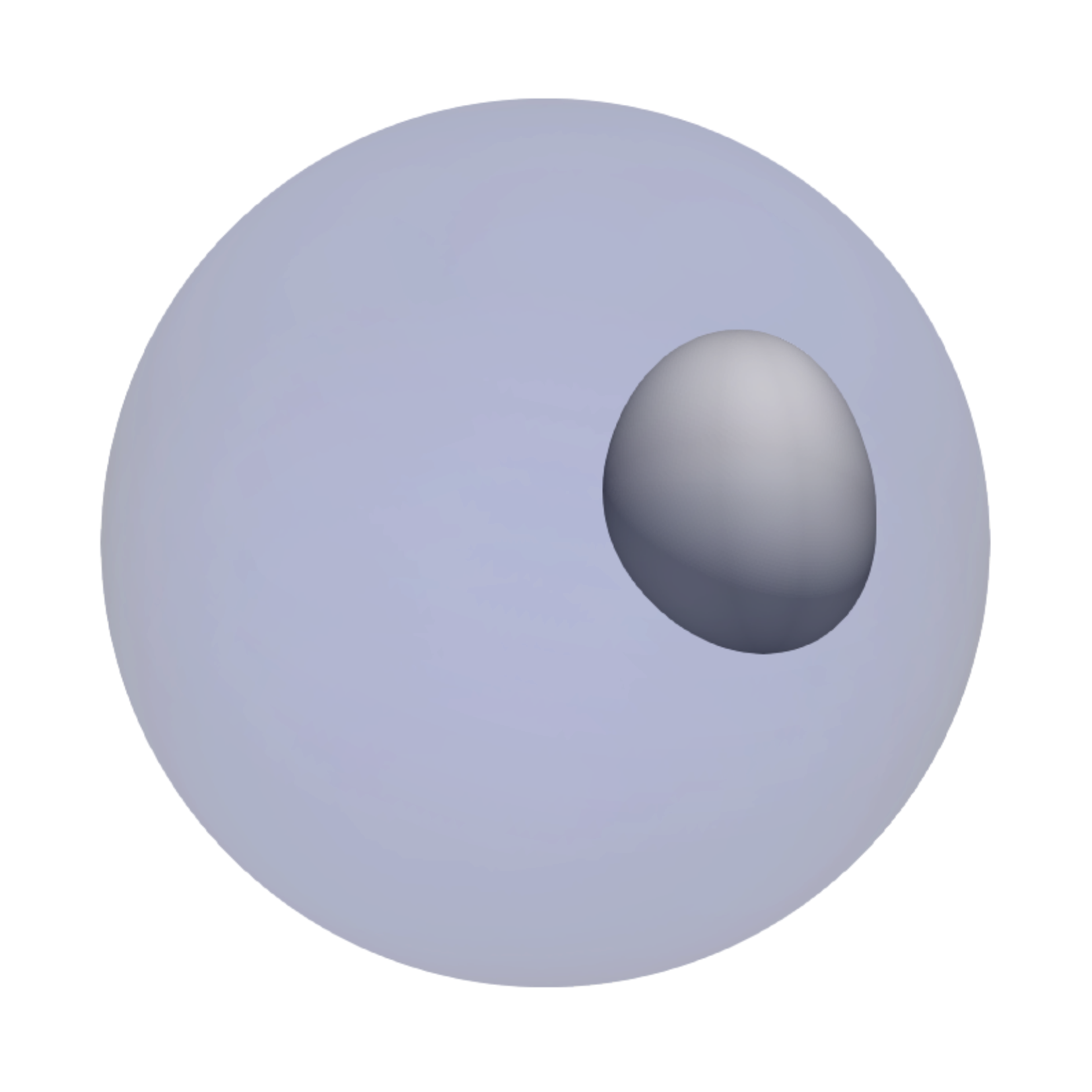}}}
\subfloat[Iteration 45]{%
\resizebox*{3.7cm}{!}{\includegraphics{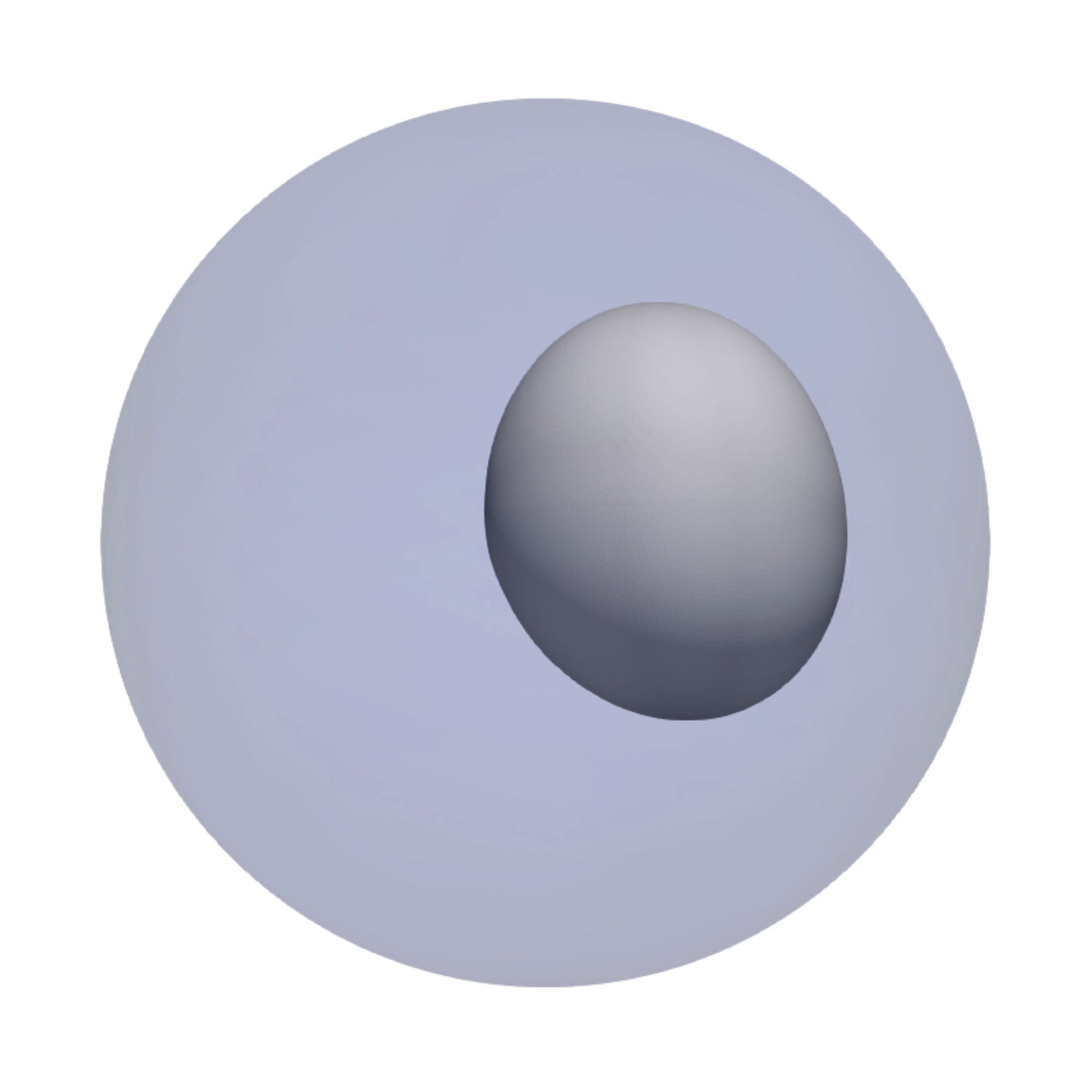}}} \vspace{3pt}
\subfloat[Iteration 60]{%
\resizebox*{3.7cm}{!}{\includegraphics{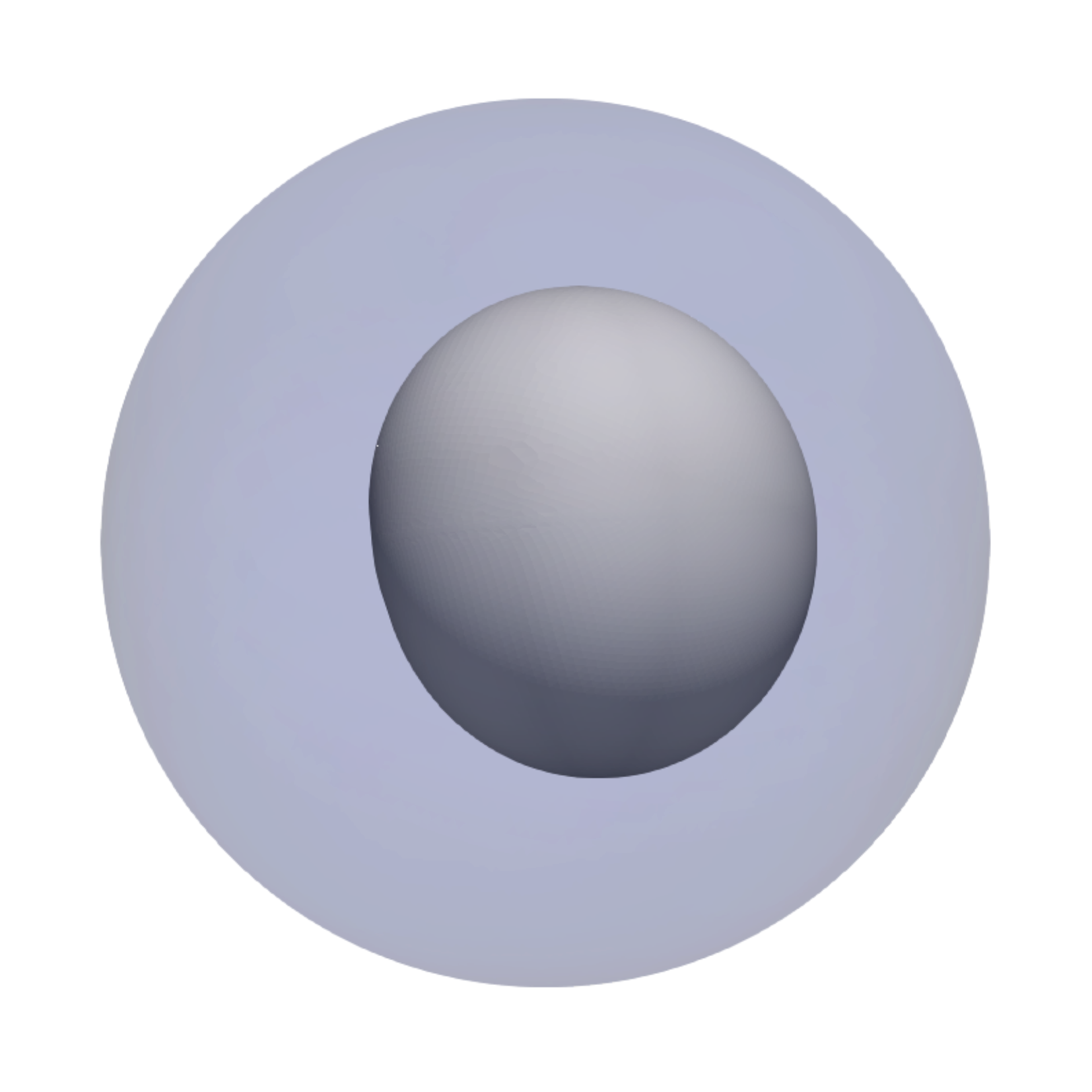}}}
\subfloat[Iteration 70]{%
\resizebox*{3.7cm}{!}{\includegraphics{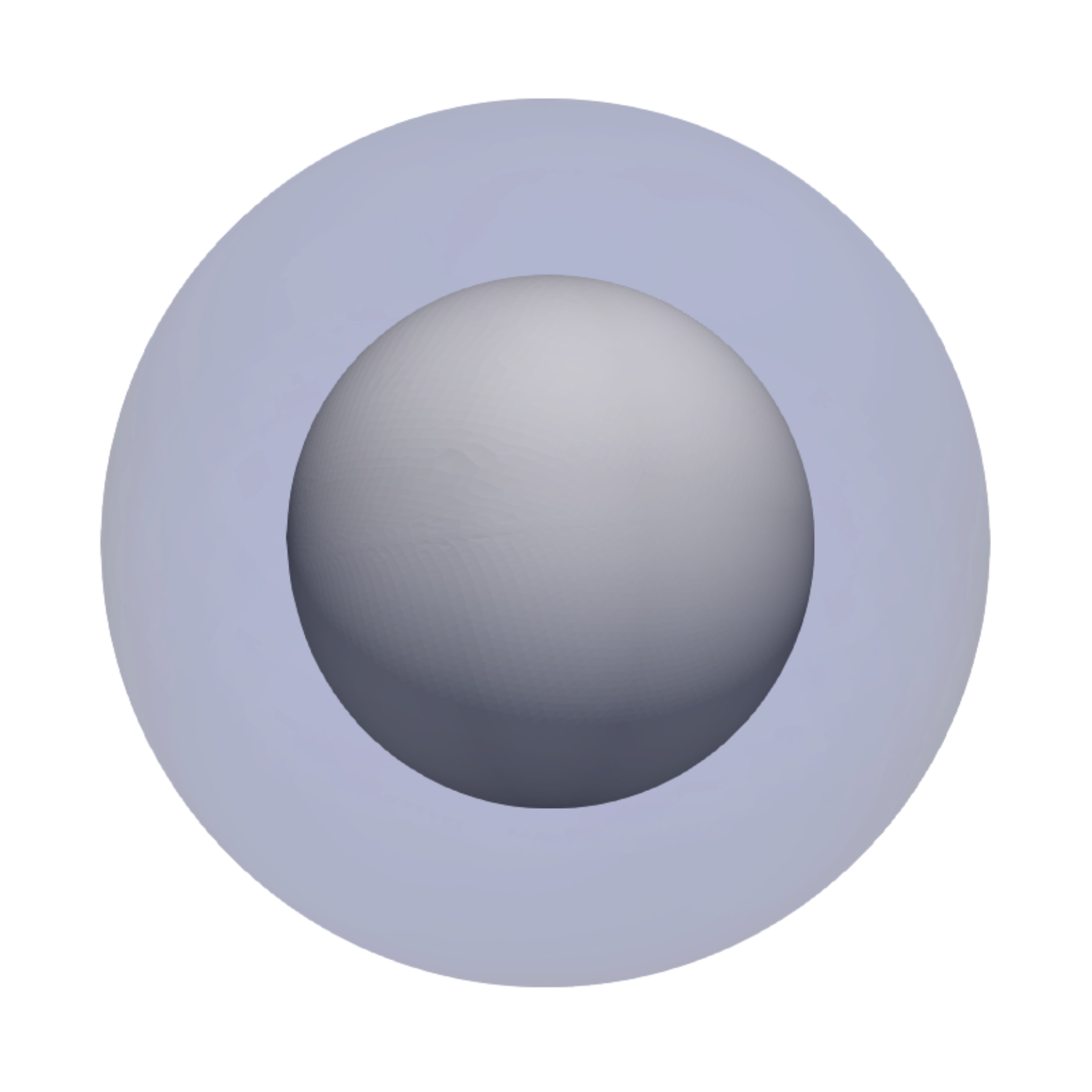}}}
\subfloat[Iteration 80]{%
\resizebox*{3.7cm}{!}{\includegraphics{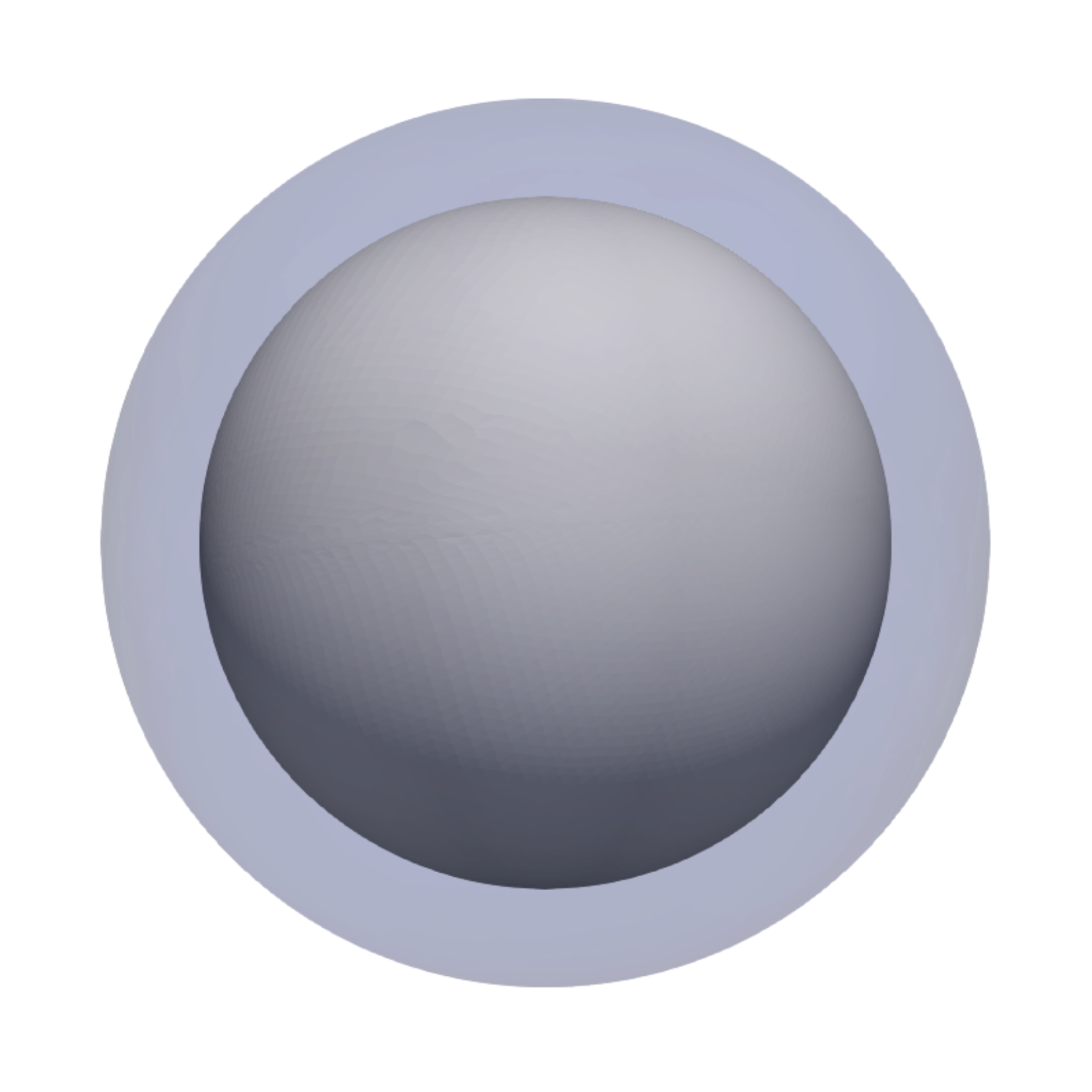}}}
\subfloat[Iteration 101]{%
\resizebox*{3.7cm}{!}{\includegraphics{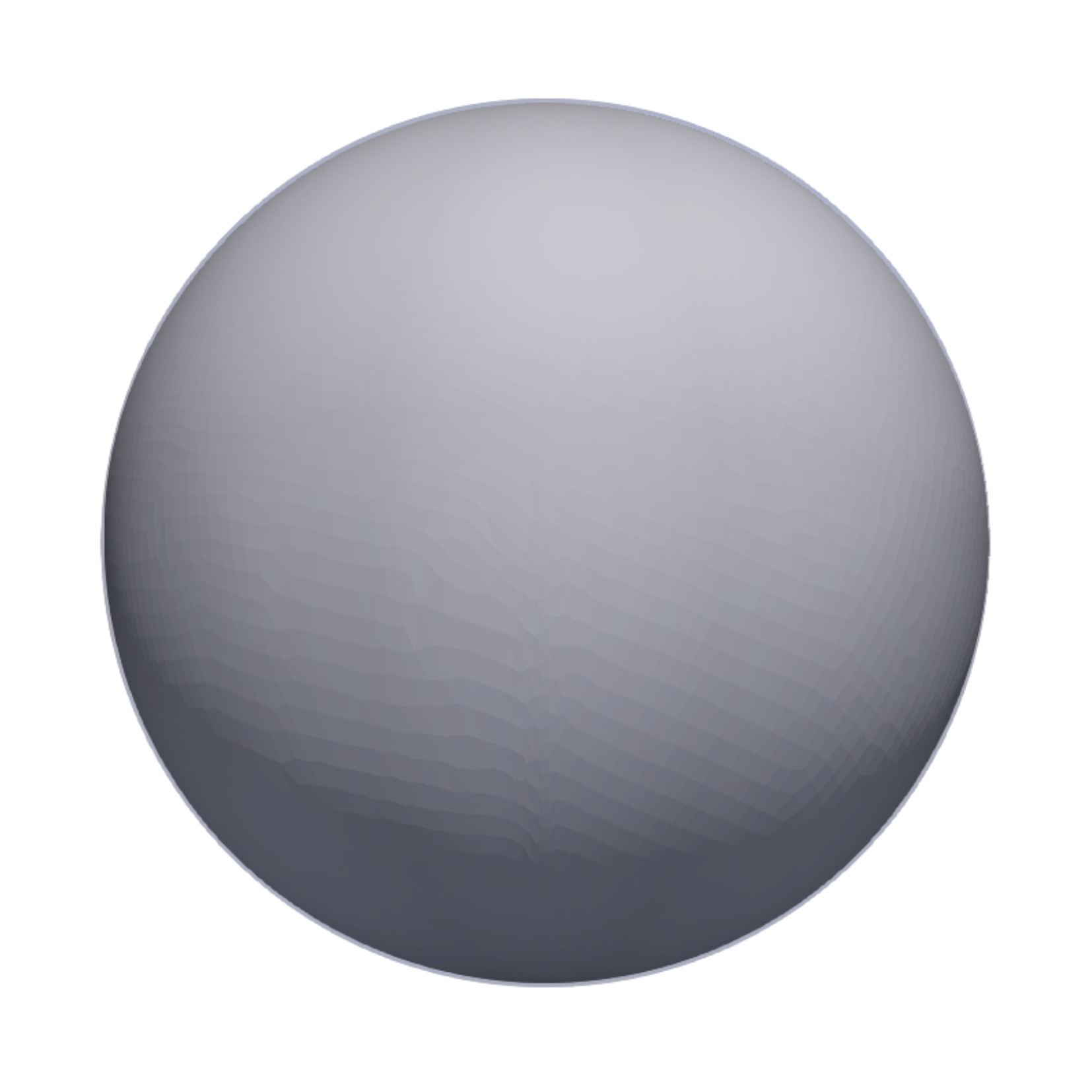}}}
\caption{Design updates.} \label{fig:sphere_iters}
\end{figure}

\subsection{Real-world example}
We consider a real-world application to shape optimization. The present method is implemented in ShapeModule, which is a flexible solver-agnostic optimization platform and provides optimization algorithms as well as shape control methods, such as Vertex Morphing \cite{Hojjat}. The optimization problem is to minimize the mass of a frame structure under load-displacement constraint (i.e., the displacement of every surface node is bounded). The optimization problem writes
\begin{equation}
\begin{split}
& \textnormal{minimize} ~~~~   M(x), \\
& \textnormal{subject to} ~~~  g_i(x,u) \leq 0, ~ i = 1,...,m,
\end{split}
\end{equation}
where $M(x)$ is the function for the mass, $g_i(x,u)$ is a point-wise formulated displacement constraint for the $i$-th node, $m$ is the number of nodes of the design surface mesh, $x \in \mathbb{R}^{3m}$ is the field of nodal coordinates of the design surface mesh, and $u \in \mathbb{R}^{3m}$ is the nodal displacement field. The number of design variables is 144423, and the number of constraints is 48141. Note that for multiple constraints, we can use the logarithmic barrier function formulation as in the previous test examples. Each single displacement constraint gradient can be efficiently computed using the adjoint sensitivity analysis. In this application, we use the load-displacement sensitivity provided by the software OptiStruct to conform with a standard industrial design chain. We choose the parameter $\zeta = 0.95$. We note that there is a rather large gap between our theoretical analysis and the large-scale implementation of the method, as our results are asymptotic and the oracle is a black-box. Nevertheless, we were able to obtain results that are qualitatively in agreement with our theory.

\begin{figure}
\centering
\subfloat[{The initial frame design}]{%
\resizebox*{7.5cm}{!}{\includegraphics{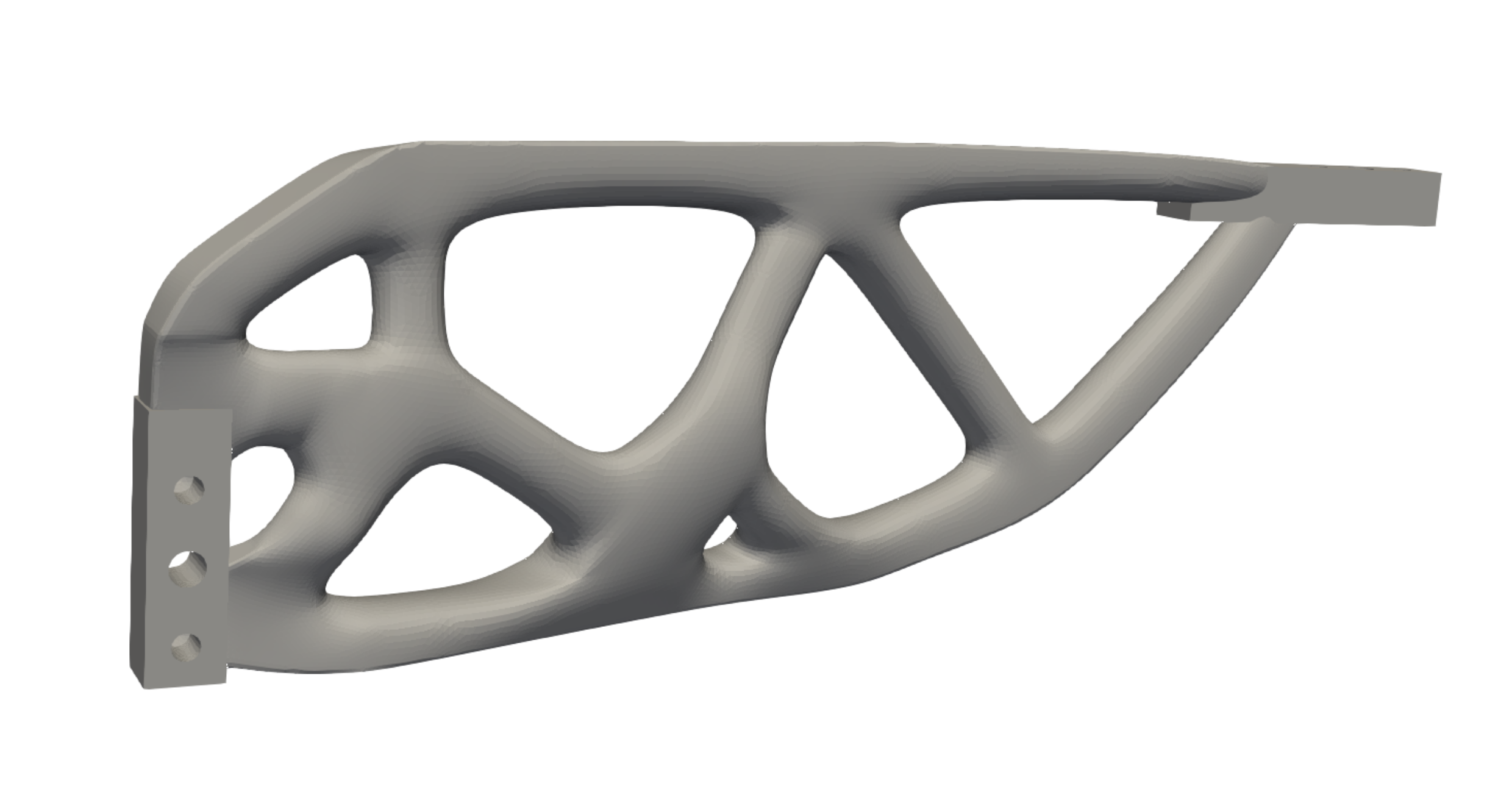}}}
\subfloat[The optimized frame design]{%
\resizebox*{7.5cm}{!}{\includegraphics{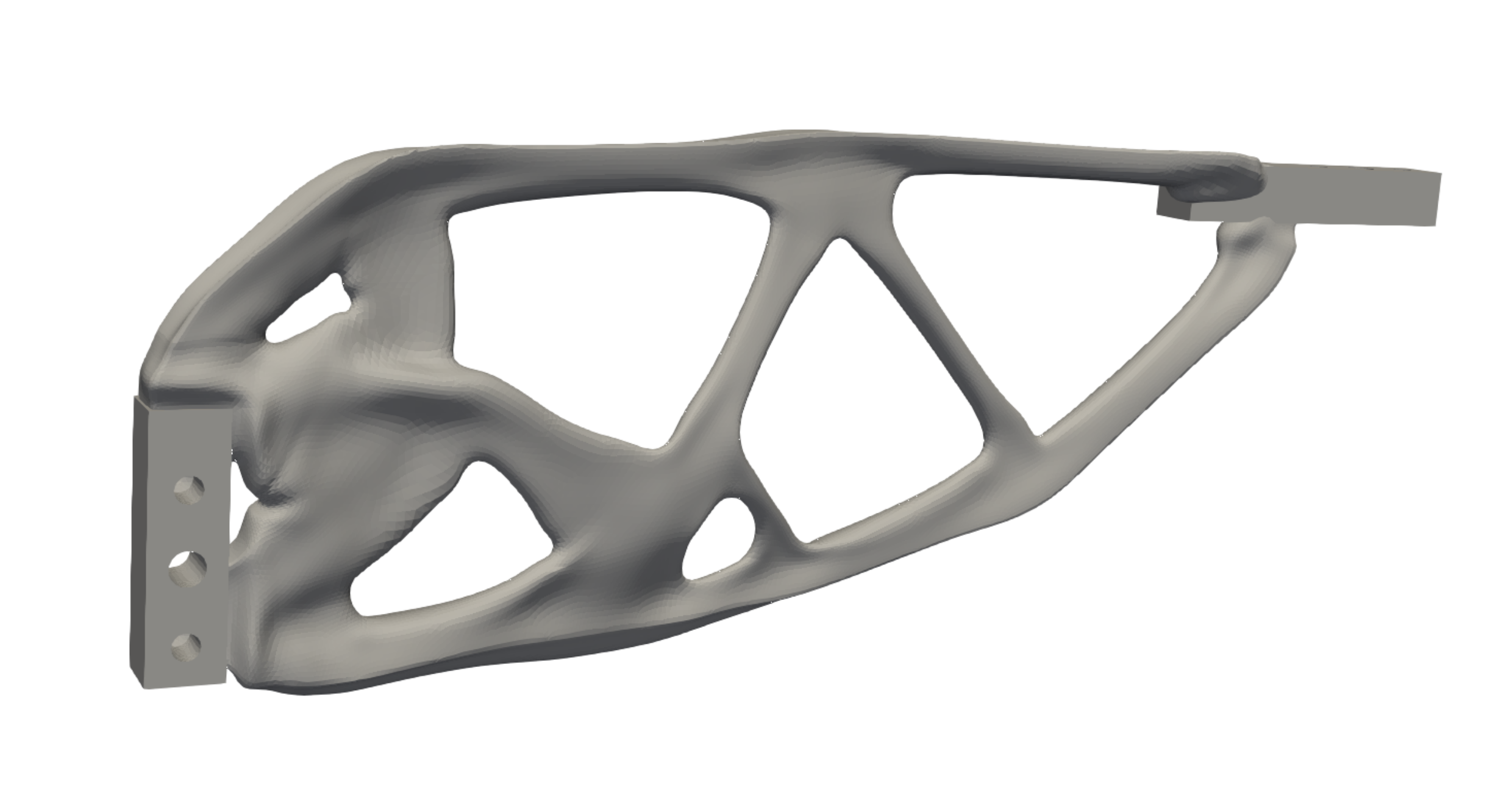}}}
\caption{Design optimization of a real-world frame structure.} \label{fig:frame}
\end{figure}

In figure \ref{fig:frame} we show the initial frame design and the optimized shape design after 194 iterations. The mass of the structure is reduced by $ 41\% $ as shown in figure \ref{fig:frame_objective}. In figure \ref{fig:frame_constraint_plot}, we show a plot of maximum constraint value $g = max\{g_i\}$ of each iteration. In figure \ref{fig:frame_centrality}, we show that the optimization is able to approach and follow a central path within the $\zeta$-neighborhood.

\begin{figure}[H]
\subfloat[Frame objective]{%
{\begin{footnotesize}%
	\executeiffilenewer{fig/frame_objective1.svg}{fig/frame_objective1.pdf}%
	{inkscape -z -C --file=fig/frame_objective1.svg %
		--export-pdf=fig/frame_objective1.pdf --export-latex}%
	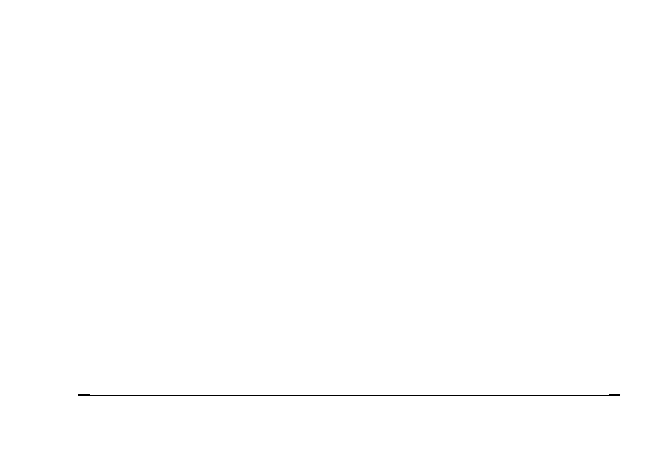%
\end{footnotesize}\label{fig:frame_objective}}} 
\subfloat[Frame constraint]{%
{\begin{footnotesize}%
	\executeiffilenewer{fig/frame_constraint_plot1.svg}{fig/frame_constraint_plot1.pdf}%
	{inkscape -z -C --file=fig/frame_constraint_plot1.svg %
		--export-pdf=fig/frame_constraint_plot1.pdf --export-latex}%
	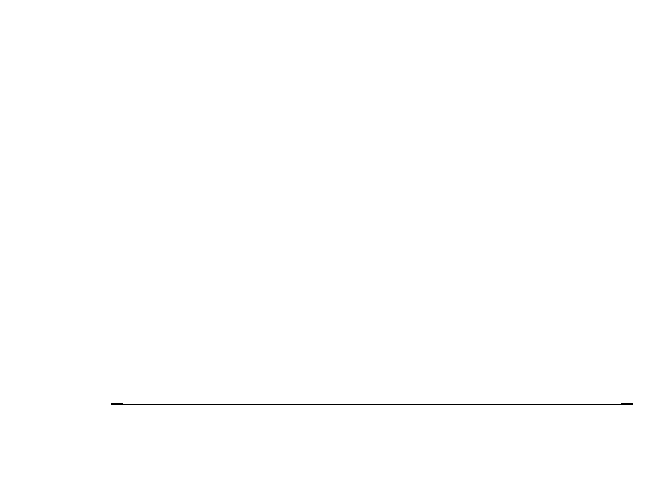%
\label{fig:frame_constraint_plot}\end{footnotesize}}}
\caption{Plot of the frame objective and constraint.} \label{fig:SCS_3000}
\end{figure}

\begin{figure}[h]
	\centering
	\begin{footnotesize}
	\executeiffilenewer{fig/frame_centrality1.svg}{fig/frame_centrality1.pdf}%
	{inkscape -z -C --file=fig/frame_centrality1.svg %
		--export-pdf=fig/frame_centrality1.pdf --export-latex}%
	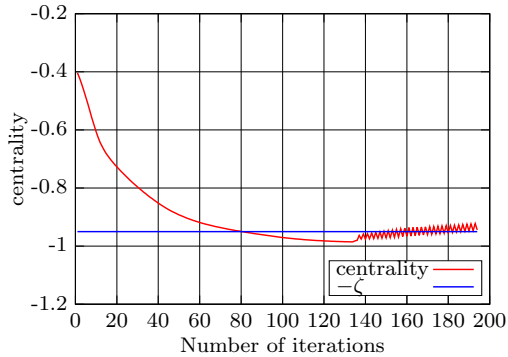%

		\caption{Plot of the centrality measure}
		\label{fig:frame_centrality}
	\end{footnotesize}
\end{figure}

\begin{remark}
By following a central path, an intermediate design improves not only the objective function but also the constraint function. Take the design of iteration 80 as an example: the mass is reduced by $20.5\%$, and the displacement is reduced by $12.1 \%$. These designs may enrich the design options if the original problem is reformulated as a bi-objective optimization problem, in which both mass and the maximum displacement are set as objectives. The resulting intermediate designs alongside a central path are approximated Pareto solutions.


\end{remark}

\vskip 2mm

\begin{remark}
This section focuses on the vanilla implementation of GDAM and demonstrates its robustness in solving node-based shape optimization problems, as robustness is considered a key to practical success for real-world shape design applications \cite{bletzinger2011free}. Accelerated GDAM can be applied to further reduce the number of iterations. However, both accelerated GDAM and Vertex Morphing introduce their own auxiliary variables into the optimization process: $\{y_k\}$ in the Algorithm \ref{alg:2} for GDAM and the control variables $\mathbf{p}$ in Vertex Morphing, respectively. Moreover, the forward and backward mapping between the control space and variable space in VM makes an application of accelerated GDAM even more delicate. The most efficient and robust way to combine both methods is worthy of further investigation. And this is of independent interest than this manuscript, which we leave to future work.
\end{remark}

\section{Application to Sensor Network Localization via Semidefinite Programming}
\label{sec:snl}
In this section, we apply the present method to solve sensor network localization (SNL) problems via the semidefinite programming (SDP) formulation \cite{biswas2004semidefinite}. First, we briefly introduce the basics of SDP-based SNL, after which we show computational experiments and comparisons with state-of-the-art first- and second-order SDP solvers.

\subsection{Sensor Network Localization via Semidefinite Relaxation}

Sensor Network Localization (SNL) aims to recover unknown sensor locations with (partially) given distance information between them, and thus can be considered an inverse problem. Consider a nonlinear least square (NLS) formulation for SNL problems with $n$ sensors and $m$ anchors in 2d:
\begin{equation}
\tag{NLS}
\min_{x_1,..., x_n \in \mathbb{R}^2} \sum_{(i,j) \in \mathcal{M} } \left( | x_i - x_j |^2   - d_{ij}^2 \right)^2 + \sum_{(k,j) \in \mathcal{\bar{M}} } \left( | a_k - x_j |^2   - d_{kj}^2 \right)^2,
\label{eq:SNL_LS}
\end{equation}
where the location of each sensor $x_i \in \mathbb{R}^d, i = 1,..., n,$ is to be determined, and the location of each anchor $a_k \in \mathbb{R}^2, k = 1,..., m$ is known. The distance $\{ d_{ij} : (i,j) \in \mathcal{M} \} $ and $\{ d_{kj} : (k,j) \in \mathcal{\bar{M}} \} $ are known observations between sensor-sensor pairs and sensor-anchor pairs, respectively. Often, the distances are detected in a given radius $r$.

We can formulate SNL as a constrained optimization problem,
\begin{equation}
\begin{split}
\min_{x_1,..., x_n \in \mathbb{R}^2} ~~ & 0, \\
\text{s.t.} ~~~ &| x_i - x_j |^2   = d_{ij}^2, ~~ (i,j) \in \mathcal{M}  \\
& | a_k - x_j |^2   = d_{kj}^2, ~~ (k,j) \in \mathcal{\bar{M}} \\
\end{split}
\label{eq:SNL_FF}
\end{equation}

SNL is a typical nonconvex quadratic optimization problem \cite{luo2010semidefinite}\cite{nesterov2000semidefinite} for which one is not content with finding a local solution but the global optimum, so that the true sensor locations can be recovered. Following the work \cite{biswas2004semidefinite}, we relax the NP-hard nonconvex problem \eqref{eq:SNL_FF} using SDP. First, we rewrite the distances of sensor-sensor pairs,
\begin{equation}
    | x_i - x_j |^2 = (e_i - e_j)^T X^T X (e_i - e_j) = \langle E_{ij},  X^T X \rangle,
\end{equation}
where $e_i \in \mathbb{R}^n$ is the $i$-th unit vector, $X$ is a $2 \times n$ matrix whose $i$-th column is $x_i$, and $
E_{ij} =  (e_i - e_j)  (e_i - e_j)^T \in \mathbb{S}^n$. Similarly, we have
\begin{equation}
| a_k - x_j |^2 = \left[ a_k^T,  -e_j^T
\right]
\begin{bmatrix}
I_2 & X\\
X^T & X^T X
\end{bmatrix} \begin{bmatrix}
a_k \\
-e_j 
\end{bmatrix}  
= \langle \bar{M}_{kj}, Z \rangle,
\end{equation}
where
\begin{equation*}
\bar{M}_{kj} = \begin{bmatrix}
a_k \\
-e_j 
\end{bmatrix} \left[ a_k^T,  -e_j^T \right],
\end{equation*}
and
\begin{equation}
Z = 
\begin{bmatrix}
I_2 & X\\
X^T & X^T X
\end{bmatrix}.
\end{equation}
Furthermore, let
\begin{equation}
M_{ij} = 
\begin{bmatrix}
0 & 0\\
0 & E_{ij}
\end{bmatrix}.
\end{equation}

By lifting the variables, we obtain an equivalent formulation in the variable $Z \in \mathbb{S}^{n+2}$ for \eqref{eq:SNL_FF},
\begin{equation}
\begin{split}
\min ~~ & 0, \\
\text{s.t.} ~~~ & \langle M_{ij}, Z \rangle = d_{ij}^2, ~~ (i,j) \in \mathcal{M}  \\
& \langle \bar{M}_{kj}, Z \rangle = d_{kj}^2, ~~ (k,j) \in \mathcal{\bar{M}} \\
& Z_{1:2,1:2} = I_2, ~Z \succeq 0,\\
& \text{rank}(Z) = 2.
\end{split}
\label{eq:snl_sdp}
\end{equation}
The problem \eqref{eq:snl_sdp} is, still, a very challenging nonconvex optimization problem. However, by lifting the variable, the nonconvexity is now only allocated at the rank constraint
\begin{equation*}
 \text{rank}(Z) = 2.
\end{equation*}
The semidefinite relaxation is to drop out this rank constraint, and thus a convex optimization problem is obtained,
\begin{equation}
\tag{SNLP}
\begin{split}
\min ~~ & 0, \\
\text{s.t.} ~~~ & \langle M_{ij}, Z \rangle = d_{ij}^2, ~~ (i,j) \in \mathcal{M}  \\
& \langle \bar{M}_{kj}, Z \rangle = d_{kj}^2, ~~ (k,j) \in \mathcal{\bar{M}} \\
& Z_{1:2,1:2} = I_2,~ Z \succeq 0.\\
\end{split}
\label{eq:snl_sdr}
\end{equation}
The dual of \eqref{eq:snl_sdr} writes,
\begin{equation}
\tag{SNLD}
\begin{split}
\max_{V,y} ~~ & \langle I_2, V \rangle + \sum_{(i,j)\in \mathcal{M}} y_{ij} d_{ij}^2 + \sum_{(k,j)\in \mathcal{\bar{M}}} y_{kj} d_{kj}^2\\
\text{s.t.} ~~~ & \begin{bmatrix}
V & \mathbf{0}\\
\mathbf{0} & \mathbf{0}
\end{bmatrix} + \sum_{(i,j)\in \mathcal{M}} y_{ij} M_{ij} + \sum_{(k,j)\in \mathcal{\bar{M}}} y_{kj} \bar{M}_{kj} + S = 0, \\
&S \succeq 0.\\
\end{split}
\label{eq:snl_sdr_dual}
\end{equation}

The semidefinite relaxation is mathematically elegant, and many applications show that, indeed, SDRs result in tight approximations to the original NP-hard problems, see, e.g., the seminal work on the maximum cut problem \cite{goemans1995improved}. In the case of SNL, the work \cite{so2007theory} shows that SDR provides the exact solution to the original problem if the sensors are uniquely localizable with provided distance information. On the other side, the cost of SDR is the much-increased problem dimensions and the introduced nontrivial PSD cone constraint, which are computationally challenging for large-scale applications.

Notice that the SDR formulation \eqref{eq:snl_sdr} of an SNL problem is a feasibility program \textemdash any feasible solution solves the optimization problem. The present method is not directly applicable since the optimization iterates in the feasible set. On the other hand, the dual formulation \eqref{eq:snl_sdr_dual} is an inequality constrained problem and is, therefore, a well-suited optimization formulation for the present method. We implement the accelerated GDAM (Algorithm \ref{alg:2}) in MATLAB R2022a to solve the dual problem \eqref{eq:snl_sdr_dual}. Further details on implementation are reported in Appendix \ref{appendix:2}. 

\subsection{Sensor Network Localization examples}
We randomly generate a set of examples with $n$ nodes representing the true sensor locations $\left[  -0.5,  0.5 \right]^2$, with $n$ ranging from 100 to 10000. Additional four anchor nodes are generated whose positions are known \textit{a priori}. 
A total of $m$ pairs of distances information are measured within a radio detection range $r$. To measure the estimation accuracy, we use the root-mean-square-distance (RMSD), which has been widely used on SNL to test the performance and is defined as
\begin{equation}
    \text{RMSD} = \left( \frac{1}{n} \sum_{i=1}^n |x_i - x_i^\text{true} |^2 \right)^{1/2},
\end{equation}
where $x_i$ is the estimated position and $x_i^\text{true}$ is the true position of sensor $i$, respectively. All experiments are run on a Linux workstation with a 2.70 GHz 6-core AMD Ryzen 5 5600U processor and 16 GB RAM.

\subsubsection{GDAM results}
\label{sec:GDAM_SNL}
We report the runtimes and RMSDs of the present method in table \ref{tb:SNL_GDAM}. It can be seen that the obtained solutions are of moderate accuracy in terms of the RMSD. The eigenvalues of the solution matrices $Z$ show that they are near rank 2 matrices. Thus, the present method not only solves the SDR problem \eqref{eq:snl_sdr} approximately, but also solves the rank-constrained SDP problem \eqref{eq:snl_sdp} to moderate accuracy. As \eqref{eq:snl_sdp} is equivalent to the original feasibility problem \eqref{eq:SNL_FF}, we can hope that the present method provides near-optimal solutions to the original NP-hard SNL problems. Indeed, as we observed in our extensive numerical experiments, using the obtained SDR solutions as initializations, we can accurately recover all the true sensor locations by applying local search methods such as gradient descent for \eqref{eq:SNL_LS}. We show some of the graphical results in figure \ref{fig:gdam_snl_3k} and \ref{fig:gdam_snl_10k}. The left figures show the results obtained by applying accelerated GDAM for the SDR \eqref{eq:snl_sdr_dual}; the right figures show refinements of the GDAM solutions by applying gradient descent to the original and low-dimensional problem \eqref{eq:SNL_LS}.

\begin{table} [h]
\begin{center}
\caption{GDAM results for SNL problems.} 
	\begin{tabular}{|c |c |c | c | c | c | c |}
		\hline
		 n & m & r &$\zeta$  & iters. & RMSD & runtime (hh:mm:ss)  \\ [0.5ex]
		\hline\hline
		100  & 1058 & 0.3 & 0.9999 & 256  & 3.09e-3 & $ <1s$  \\  
        \hline
		500& 14699 & 0.21 & 0.9999 & 948 & 2.29e-3 & 00:00:09  \\  
		\hline
		1500  & 98007 & 0.18 & 0.9999 & 1312 & 3.27e-3 & 00:02:41   \\  
		\hline  
		3000 & 246177 & 0.14 & 0.9999 & 1164   & 5.38e-3 & 00:12:37    \\  
        \hline
		5000  & 510116 & 0.12 & 0.9999 & 1231  & 5.63e-3 & 00:49:58     \\  
		\hline 
		10000 & 1446144 & 0.10 & 0.9999 & 1427   & 3.97e-3 & 05:51:51     \\  
		\hline    
	\end{tabular}\label{tb:SNL_GDAM}
 \end{center}
\end{table}

\begin{figure}[H]
\raggedleft
\subfloat[GDAM for SDR]{%
\resizebox*{7.5cm}{!}{\includegraphics{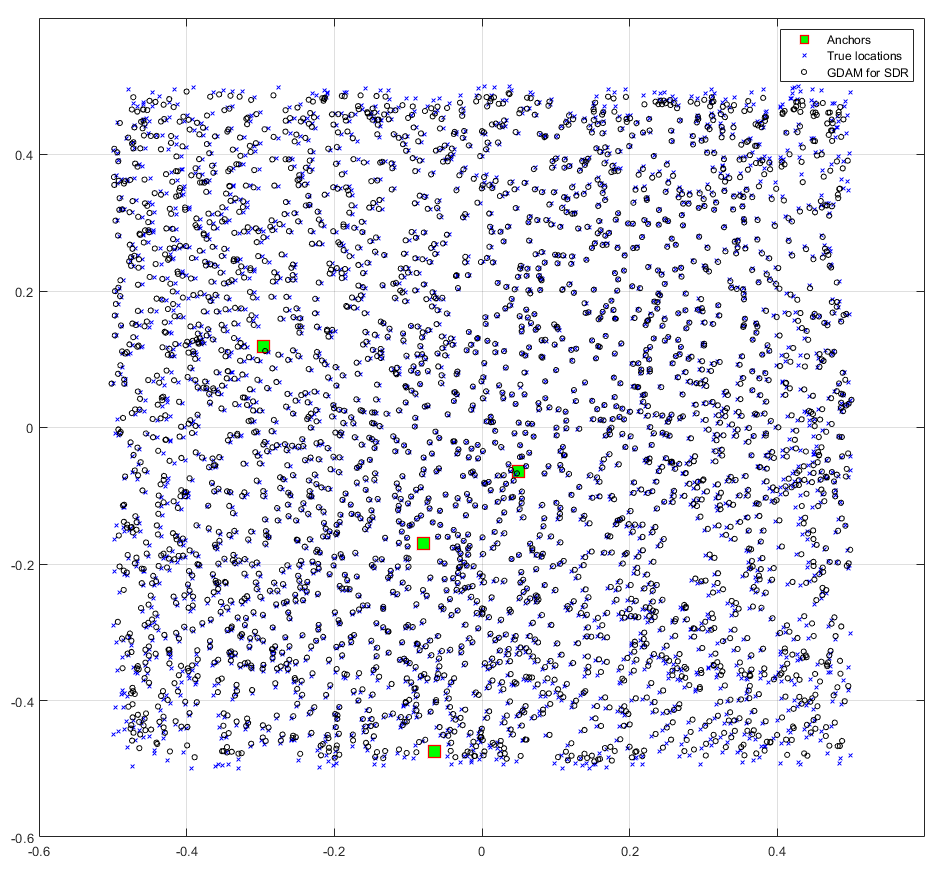}}} 
\subfloat[GD after SDR]{%
\resizebox*{7.5cm}{!}{\includegraphics{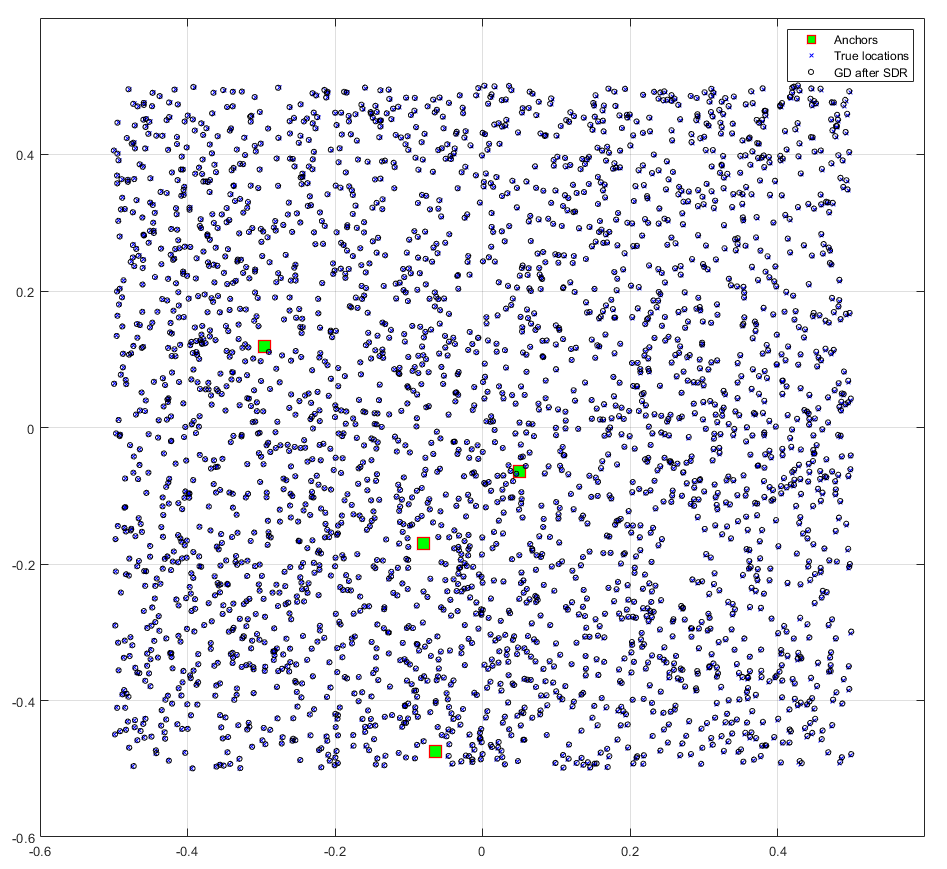}}}
\caption{Graphical results of GDAM for the SNL problem with 3000 sensors.} \label{fig:gdam_snl_3k}
\end{figure}

\begin{figure}[h]
\subfloat[GDAM for SDR]{%
\resizebox*{7.5cm}{!}{\includegraphics{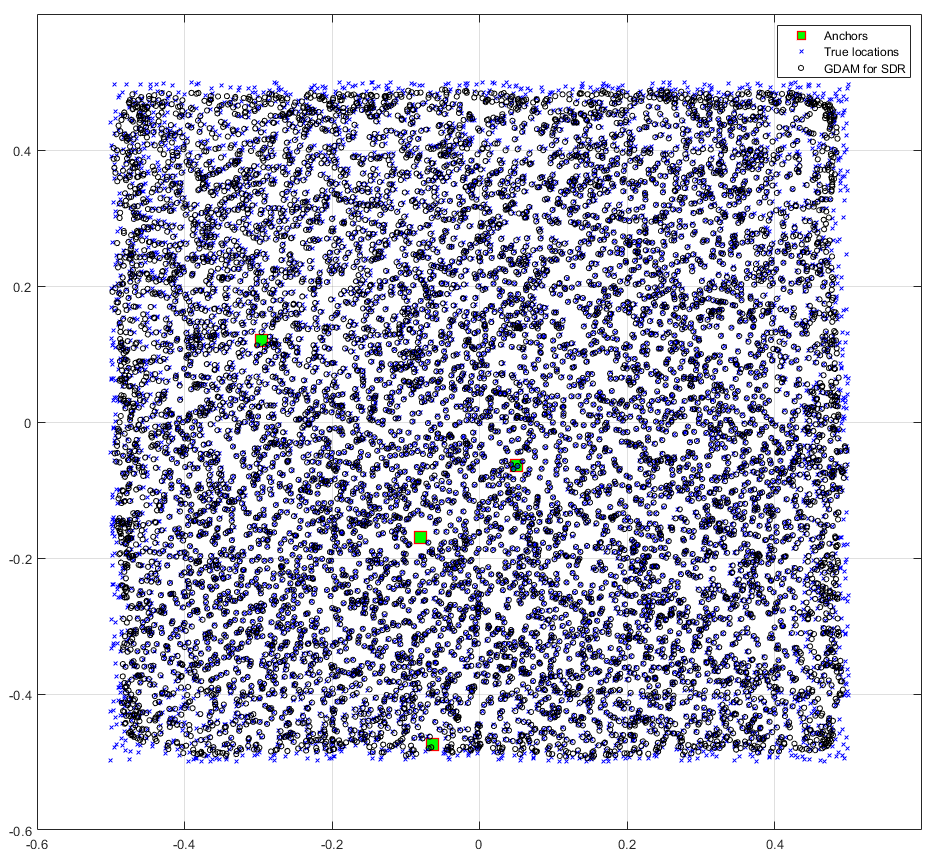}}}
\subfloat[GD after SDR]{%
\resizebox*{7.5cm}{!}{\includegraphics{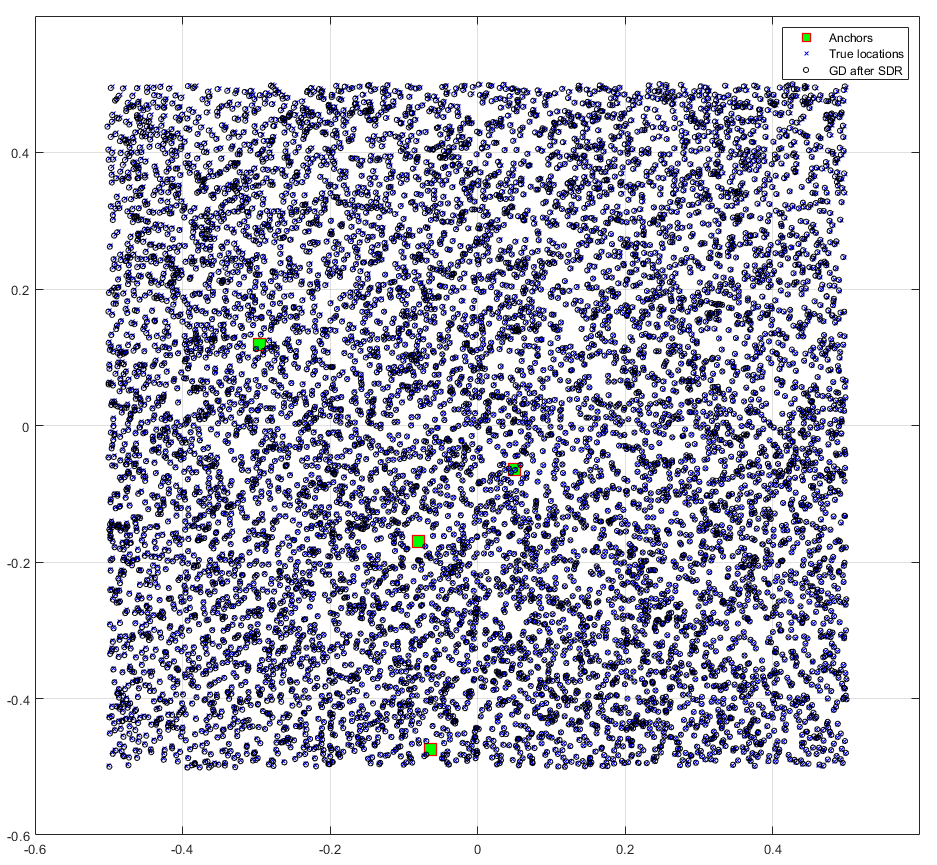}}}
\caption{Graphical results of GDAM for the SNL problem with 10000 sensors.} \label{fig:gdam_snl_10k}
\end{figure}

\subsubsection{Comparisons}
 We compare results in terms of runtime and RMSD with three well-established solvers: DSDP v5.8 \cite{benson2008algorithm}, a generic SDP solver based on a second-order dual interior-point method,  SCS (Splitting Conic Solver) version 3.1.0 \cite{odonoghue21}, a generic large-scale conic solver based on ADMM, and SDPNAL+ version 1.0 \cite{sun2020sdpnal+}, a large-scale SDP solver based on a semismooth Newton-CG augmented Lagrangian method.  We run the same examples listed in table \ref{tb:SNL_GDAM}.

\vskip 2mm

\noindent \textbf{DSDP}:
We use the default settings of DSDP, where the tolerance for the optimality gap is set to be $1e-6$. DSDP is an efficient second-order solver for SNL problems that exploits the problem feature which leads to reduced storage and increased efficiency. Table \ref{tb:sdp_dsdp} shows that DSDP can solve problems to very high accuracy. For problems with size $n \geq 1500$, the workstation is running out of memory. While DSDP achieves high accuracy for the SNL problems, GDAM requires a much less time in finding an approximate solution. To recover the true sensor locations, we run GD for the low-dimensional nonlinear least squares problem \eqref{eq:SNL_LS}, a well-established approach for SNL problems \cite{biswas2006semidefinite} that belongs to the general two-phase strategy for solving SDP relaxation problems, where a SDP solutions is used as an initialization for nonlinear programming methods to locally solve the original NP-hard problem \cite{luo2010semidefinite}.
 
\begin{table} [h]
	\begin{center}
\caption{DSDP results for SNL problems} 
	\begin{tabular}{|c | c | c | c | c |}
		\hline
		 n & $tol.$  & iters. & RMSD & runtime (hh:mm:ss)  \\ [0.5ex]
		\hline\hline
		100  &  1e-6 & 36  & 5.06e-10 & $< 1 s$  \\  
        \hline
		500  &  1e-6 & 32 & 6.37e-10 & 00:06:47  \\  
		\hline
	\end{tabular}
\label{tb:sdp_dsdp}
\end{center}
\end{table}

\vskip 2mm

\noindent \textbf{SCS}:
We use `sparse-indirect' for the linear solve setting and choose convergence tolerances between $tol. = 1e-2$ and $tol. = 1e-4$, and set the maximum runtime to be 24 hours. We observe that SCS encounters numerical difficulties if the objective is set to be empty, as is the case of our feasibility formulation \eqref{eq:snl_sdr}. Therefore, We add an objective $ \pm\text{trace}(Z)$ to the model. We find that the runtime is quite sensible with the choice of the sign. Therefore, we tuned the sign $\pm$ and report results with better performance. We use Yalmip (Version 31-March-2021) \cite{lofberg2004yalmip} for the problem modeling and the interface to SCS.

We report the results of SCS in table \ref{tb:SNL_SCS}. For problems with $n \leq 1500$, a tolerance with $1e-3$ yields comparably accurate results as GDAM. For the problem with $n = 3000$, a smaller tolerance $1e-4$ is needed to obtain a good estimation as initialization for local refinement, as shown in figure \ref{fig:SCS_3000}. For all the instances, GDAM is more competitive in the runtime. For example, for the problem with $n = 3000$, GDAM takes about 12 minutes to find a moderate accurate solution with RMSD $= 5.38 e-3$, and SCS finds an inferior solution after five and a half hours.

\begin{table}[h]
 \begin{center}
 \caption{SCS results for SNL problems.}
\begin{tabular}{|c|c|c|c|c|}
\hline
n & $tol.$ & iters. & RMSD & runtime (hh:mm:ss)  \\ \hline\hline
\multirow{ 2}{*}{100} & 1e-2 & 400 & 7.64e-2 & 00:00:02 \\
  & 1e-3 & 650 & 2.91e-3 & 00:00:03 \\ \hline
\multirow{ 2}{*}{500} & 1e-2 & 400 & 1.57e-01 & 00:00:48 \\
  & 1e-3 & 950 & 1.03e-02 & 00:01:23 \\ \hline
\multirow{ 2}{*}{1500} & 1e-2 & 400 & 2.29e-1 & 00:13:21
 \\
& 1e-3 & 1225 & 1.82e-03 & 00:40:50 \\ \hline
\multirow{ 2}{*}{3000} & 1e-3 & 1950 & 1.21e-1 &  05:28:20
 \\
& 1e-4 & 7225 & 6.14e-4 & 18:35:00 \\ \hline
\end{tabular} 
\label{tb:SNL_SCS}
\end{center}
\end{table}

\begin{figure}[H]
\subfloat[SCS with $tol. = 1e-3$ for SDR]{%
\resizebox*{7.5cm}{!}{\includegraphics{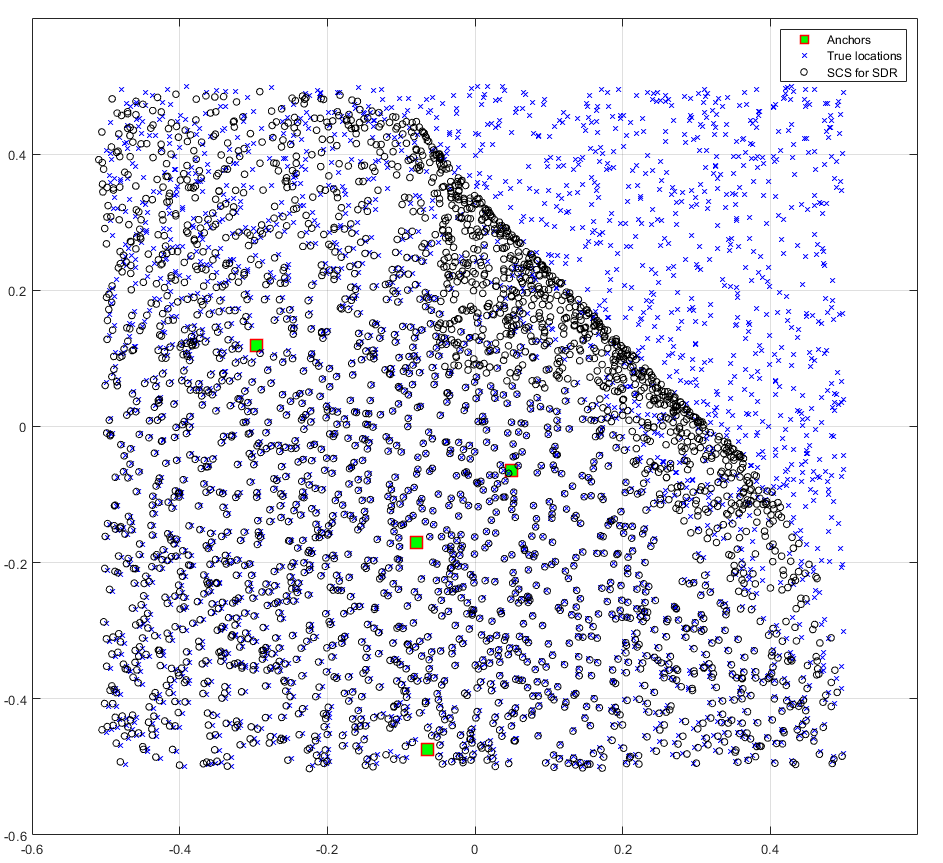}}} 
\subfloat[SCS with $tol. = 1e-4$ for SDR]{%
\resizebox*{7.5cm}{!}{\includegraphics{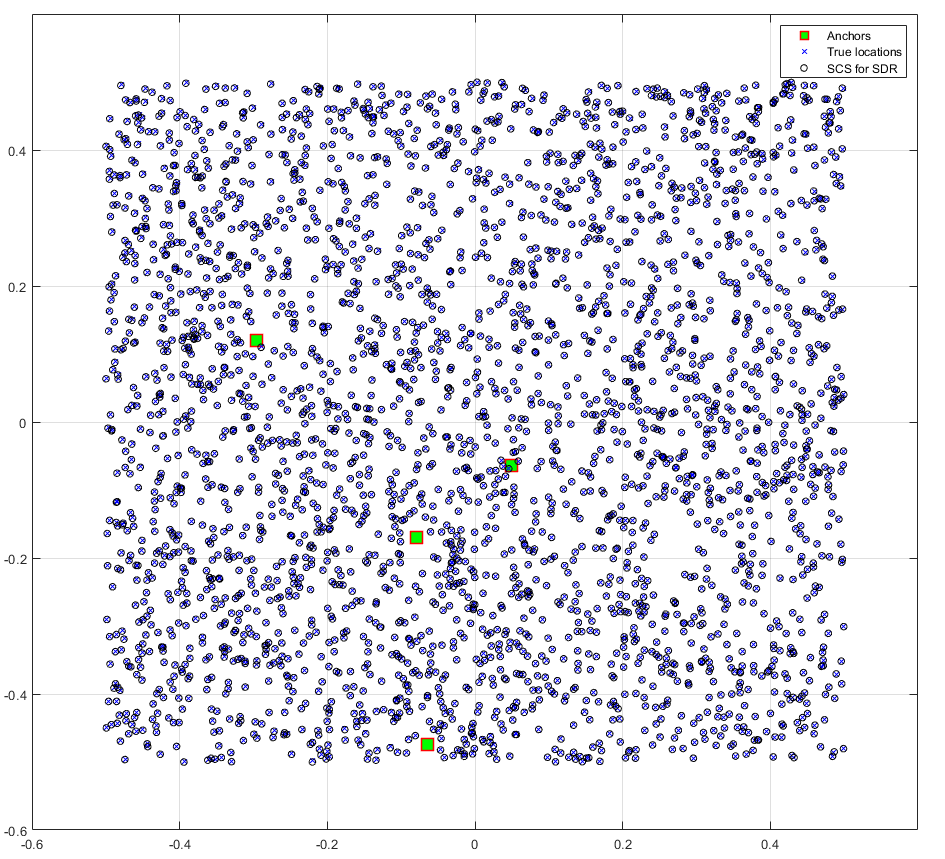}}}
\caption{Graphical results of SCS for the SNL problem with 3000 sensors.} \label{fig:SCS_3000}
\end{figure}

\vskip 2mm
\noindent \textbf{SDPNAL+}: SDPNAL+ v1.0 is a two-phase solver in which an ADMM algorithm is used in phase I to provide an initial point for the phase II algorithm, which is a semismooth Newton-CG augmented Lagrangian method. We use the default settings of SDPNAL+ and set the tolerances to $1e-2$ and $1e-3$, and set the maximum runtime to be 24 hours. We solve the test examples shown in section \ref{sec:GDAM_SNL} and report results in table \ref{tb:SNL_SDPNAL}. Comparing table \ref{tb:SNL_SDPNAL} and table \ref{tb:SNL_SCS}, we see that the performance of SDPNAL+ and SCS are fairly close. In about the same runtime, SDPNAL+ finds more accurate solutions for problems with $n\leq 1500$, while SCS finds a more accurate solution for the problem with $n = 3000$ (see SDPNAL+ result in figure \ref{fig:sdpnal_3000}). A comparison of table \ref{tb:SNL_SDPNAL} and table \ref{tb:SNL_GDAM} shows that GDAM is more efficient than SDPNAL+ in finding moderate accurate solutions in terms of RMSD for the considered SNL problems.

\begin{remark}
Large-scale semidefinite programs arising from the relaxation of SNL problems are challenging convex optimization problems. Even advanced large-scale SDP solvers, such as SCS and SDPNAL+, need to make a compromise between the computational time and the accuracy of a solution. A practical difficulty lies in how to choose such a tolerance \textit{a priori}. Take the example of $n= 3000$, SCS requires $tol. = 1e-4$ to obtain a moderate accurate localization solution, which is an order of magnitude smaller than that required for problems with $n \leq 1500$. For the same problem $n = 3000$, a practical tolerance for SDPNAL+ may lie between $1e-2$ and $1e-3$, whereas the difference in the runtime of the two tolerances is about 18 hours. On the other hand, without the knowledge of the true sensor locations, it may not be straightforward to design a practical stopping criterion. In this regard, GDAM shows its practical advantage by only using a single fixed parameter $\zeta = 0.9999$ (universally for different numbers of sensors, radius, and anchor locations). In theory, $\zeta$ is related to solution accuracy (see Theorem \ref{theorem:global_error_bound}), but its practical implications for large-scale SDP problems may be more profound, which is worth further investigation.
\end{remark}


\begin{remark}
Compared to SCS and SDPNAL+, the number of iterations needed for GDAM increases only slowly as the dimension of the problem grows. Indeed, GDAM is the only method we are aware of that is capable of solving the SDP model \eqref{eq:snl_sdr} (or \eqref{eq:snl_sdr_dual}) to 10000 sensors on a laptop workstation. This clearly shows the potential of GDAM as a practical optimization method for some of the large-scale and difficult constrained optimization problems. 
\end{remark}

\begin{table}[h]
\begin{center}
\caption{SDPNAL+ results for SNL problems.}    
\begin{tabular}{|c|c|c|c|c|}
\hline
n & $tol.$ & it. ADMM+/SSN & RMSD & runtime (hh:mm:ss)  \\ \hline\hline
100 & 1e-2 & 114/9 & 8.013e-4 & 00:00:01 \\ \hline
\multirow{ 2}{*}{500} & 1e-2 & 370/21 & 8.38e-3 & 00:00:42  \\
  & 1e-3 & 500/22 & 1.43e-3 & 00:00:47  \\ \hline
\multirow{ 2}{*}{1500} & 1e-2 & 300/31 & 4.90e-2 & 00:20:43
 \\
& 1e-3 & 3600/43 & 1.85e-5 & 00:36:53 \\ \hline
\multirow{ 2}{*}{3000} & 1e-2 & 414/38 & 2.15e-1 &  01:55:06  \\
& 1e-3 & 21603/165 & 1.40e-3& 19:49:07  \\ \hline
\end{tabular}
\label{tb:SNL_SDPNAL}	
\end{center}
\end{table}

\begin{figure}[H]
\subfloat[SDPNAL+ with $tol. = 1e-2$ for SDR]{%
\resizebox*{7.5cm}{!}{\includegraphics{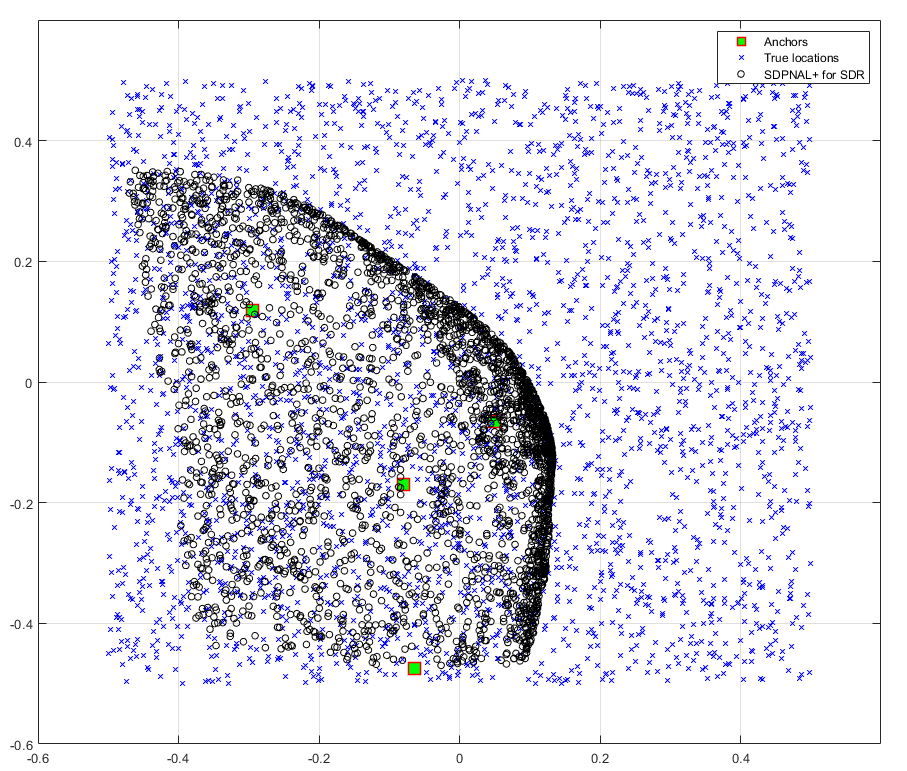}}} 
\subfloat[SDPNAL+ with $tol. = 1e-3$ for SDR]{%
\resizebox*{7.5cm}{!}{\includegraphics{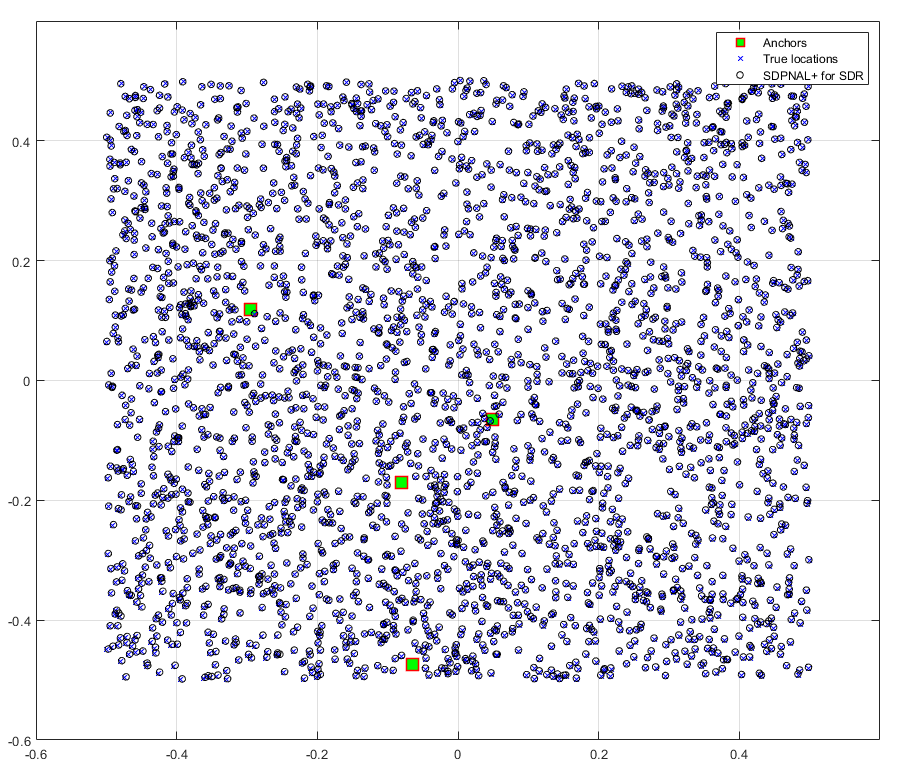}}}
\caption{Graphical results of SDPNAL+ for the SNL problem with 3000 sensors.} \label{fig:sdpnal_3000}
\end{figure}

\section{Discussion}
\label{sec:conclusion}
This paper proposes a gradient descent akin method for solving constrained optimization problems. GDAM may be considered a gradient descent method for constrained optimization and a first-order interior-point method. We give essential theoretical guarantees on the global convergence of the method to first-order solutions. We present computational algorithms based on GDAM and show their applications to two engineering problems, shape optimization and sensor network localization. Computational experiments demonstrate that GDAM is robust and very competitive in finding moderate accurate solutions and scales well to very large problems. Given the rather big difference in the characteristics of the two applications, we believe that GDAM could be suitable for other large-scale constrained optimization problems by providing inexact but useful solutions. 

\section*{Acknowledgements}

The author LC thanks the financial support of the German National High Performance Computing (NHR) Alliance. We are grateful to the ShapeModule team at the BMW Group for providing a real-world model and their shape optimization framework.

\bibliographystyle{tfs}
\bibliography{gdam_computation}

\appendix

\section{Derivation of the analytic trajectory $\Gamma^\zeta$ for problem \eqref{eq:example_proof}}
\label{appendix:1}

For problem \eqref{eq:example_proof}, the search direction field $\textbf{s}_\zeta$ reads
\begin{equation}
\begin{split}
\textbf{s}_\zeta &= - \frac{\nabla f}{|\nabla f|} - \zeta \frac{\nabla g}{|\nabla g|} \\
&= \frac{-1}{\sqrt{x_1^2 + x_2^2}}\left\lbrace  x_1,  x_2 - \zeta \sqrt{x_1^2 + x_2^2} \right\rbrace^T.
\end{split}
\end{equation}
By the uniqueness theorem of the ordinary differential equation, the integral curves of the vector field $\textbf{s}_\zeta$ has non-vanishing component $x_1$ when the initial design satisfies $x_1^0\neq 0.$ Hence we have,
\begin{equation}
\frac{dx_2}{dx_1} = \frac{x_2}{x_1} - \zeta \frac{\sqrt{x_1^2 + x_2^2}}{x_1}
\end{equation}
While $ x_1 > 0$, let $x_2 = x_1 u$, then we have $\frac{d x_2}{d x_1} = u + x_1 \frac{d u}{d x_1}$. Thus,
\begin{equation}
u + x_1 \frac{d u}{d x_1} = u - \zeta \sqrt{1 + u^2}.
\label{eq: u}
\end{equation}
We rewrite (\ref{eq: u}) as
\begin{equation}
\frac{d u}{\sqrt{1 + u^2}} = \frac{- \zeta}{x_1} d x_1.
\label{eq: du}
\end{equation}
Solving equation (\ref{eq: du}) we get
\begin{equation}
\log \left(u + \sqrt{1 + u^2}\right) = - \zeta \log x_1 + a_0.
\end{equation}
Or,
\begin{equation}
x_2 + \sqrt{x_1^2 + x_2^2} = a x_1^{1- \zeta}.
\end{equation}
Similarly, for $x_1<0$, the trajectory has the form
\begin{equation}
x_2 + \sqrt{x_1^2 + x_2^2} = a |x_1|^{1- \zeta}.
\end{equation}
Assume an initial design $\textbf{x}^0 = \left(x_1^0, x_2^0\right)$, then we get
\begin{equation}
a_\zeta = \frac{x_2^0 + \sqrt{(x_1^0)^2 + (x_2^0)^2}}{|x_1^0|^{1-\zeta}}.
\end{equation}
Hence,
\begin{equation}
x_2 + \sqrt{x_1^2 + x_2^2} = \left|\frac{x_1}{x_1^0}\right|^{1- \zeta} \left(x_2^0 + \sqrt{(x_1^0)^2 + (x_2^0)^2}\right).
\end{equation}
Let $\bar{x}_2 = \frac{1}{2} \left(x_2^0 + \sqrt{(x_1^0)^2 + (x_2^0)^2}\right)$, then we obtain the trajectory $\Gamma^\zeta$ of $\textbf{s}_\zeta$
\begin{equation}
x_2 + \sqrt{x_1^2 + x_2^2} = 2 \bar{x}_2 \left|\frac{x_1}{x_1^0}\right|^{1- \zeta}.
\end{equation}
Let  $(x_{1,\zeta}, x_{2,\zeta})$ be the point in the trajectory $\Gamma^\zeta$ with maximal $x_2$ component, then
\begin{equation}
\left\{
\begin{split}
&x_{2,\zeta} + \sqrt{x_{1,\zeta}^2 + x_{2,\zeta}^2} = 2 \bar{x}_2 \left|\frac{x_{1,\zeta}}{x_1^0}\right|^{1- \zeta},\\
&x_{2,\zeta} - \zeta\cdot \sqrt{x_{1,\zeta}^2 + x_{2,\zeta}^2} = 0.
\end{split}
\right.
\end{equation}
This gives
\begin{equation}
(x_{2,\zeta})^\zeta=\frac{2 \zeta}{1+\zeta}\frac{\bar{x}_2}{|x_1^0|^{1-\zeta}}
\left(\frac{\sqrt{1-\zeta^2}}{\zeta}\right)^{1-\zeta},
\end{equation}
and
\begin{equation}
|x_{1,\zeta}|=\frac{\bar{x}_2}{\zeta}\sqrt{1-\zeta^2}. 
\end{equation}
Let $\zeta \rightarrow 1^{-} $, then, $x_{1,\zeta} \rightarrow 0$ and $x_{2,\zeta} \rightarrow \bar{x}_2$. The trajectory $\Gamma^\zeta$ converges to the curve $\Gamma$ that is a union of the parabola
\begin{equation}
x_1^2 =  4 \bar{x}_2^2 - 4 x_2 \bar{x}_2, ~~~~ x_1\in (0, x_1^0) (~~~\mbox{or} (x_1^0,0),)
\end{equation}
and the interval $(0,\bar{x}_2)$ on $x_2$-axis.

\section{Implementation detail of GDAM for SNL}
\label{appendix:2}
We report more implementation details in addition to the section \ref{sec:implementation} for solving the SNL problems. 

\subsection{Logarithmic barrier function for the positive semidefinite cone constraint}
In the following, we briefly show how the PSD constraint can be tackled using the standard logarithmic barrier function. For the primal SDP formulation \eqref{eq:snl_sdr}, the PSD cone constraint of the matrix variable $Z \succeq 0 $ can be expressed as $Z$ only has non-negative eigenvalues, i.e., 
\begin{equation}
\lambda_i \geq 0, i = 1,...,n,
\label{eq:eigenvalue_cons}
\end{equation}
where we assume $Z \in \mathbb{S}^{n}$ for a general presentation. Using the logarithmic barrier function for \eqref{eq:eigenvalue_cons}, we have
\begin{equation}
\Phi(Z) = \sum_{i=1}^n - \log \lambda_i = -\log \left(  \Pi_{i=1}^n \lambda_i \right)= - \log \left( \det(Z) \right).
\label{eq:log_det_X}
\end{equation}
The gradient of $\Phi(Z)$ with respect to the matrix variable $Z$ writes
\begin{equation}
    \frac{\partial}{\partial Z_{ij}} \left( -\log \det Z \right) = -\left( Z^{-1} \right)_{ji}.
\label{eq:grad_log_det_X}
\end{equation}

For the dual SDP formulation \eqref{eq:snl_sdr_dual}, the PSD cone constraint can be formulated similarly. Notice now that the variable is the dual variable $\mathbf{y}$, and so we have
\begin{equation}
    \Phi(\mathbf{y}) = - \log \left[ \det(S(\mathbf{y})) \right].
\end{equation}

The gradient of $\Phi(\mathbf{y})$ writes
\begin{equation}
    \frac{\partial}{\partial y_{i}} \bigl(- \log \left[ \det(S) \right] \bigr) = \langle  S^{-1},\mathbf{A}_i \rangle_{\mathbb{S}^n},
\label{eq:grad_log_det_S}
\end{equation}
where $A = \{M_{ij}, \Bar{M}_{kl}, (i,j) \in \mathcal{M}, (k,l)\in \mathcal{\bar{M}} \}$. From \eqref{eq:grad_log_det_X} and \eqref{eq:grad_log_det_S}, we see that in order to use a gradient-based method with the logarithmic barrier framework, a matrix inversion needs to be done for a matrix in $\mathbb{S}^{n}$ for both the primal and dual formulation.

\subsection{Presolve} 
Presolve is an important part of the practical implementation of optimization algorithms that transforms the input problem into an easier one. For SNL problems, we use a simple presolve strategy that rescales the network geometry to the range $\left[-5,  5 \right]^2$. Numerically, it can be regarded as a preconditioning of the formulated dual SDP problem \eqref{eq:snl_sdr_dual}. Compared to the canonical sensor range $\left[-0.5,  0.5 \right]^2$, which is often presented the literature, the dual IPM solver DSDP justifies the efficacy of our presolve strategy with a reduced number of iterations. We note that general purpose presolve/preconditioning strategies (see, e.g., \cite{zheng2021chordal}) will probably further improve the performance of the method, which we leave for future work.

\subsection{Strictly feasible initialization}
Different initialization strategies are available for IPMs in solving SDPs. In this work, we use the Phase I-then-Phase II method, i.e., we first try to find a feasible point (and in our case one for the dual problem), and then start the main solve routine. Notice that the dual problem \eqref{eq:snl_sdr_dual} is feasible when $V = 0$ and $y_{ij} = 0$ for all $(i,j) \in \mathcal{M} $ and $y_{kj} = 0$ for all $(i,j) \in \mathcal{\bar{M}} $, and we denote this zero initialization as $y_0$ for convenience. The point $y_0$ is, however, not an interior point in the feasible set for which we can start GDAM. To overcome this, we propose to run a few steps of gradient descent, starting from $y_0$, for the following auxiliary logarithmic barrier problem,
\begin{equation}
\min ~~ - \log \det \Big( \lambda I_{n+2} - \begin{bmatrix}
V & \mathbf{0}\\
\mathbf{0} & \mathbf{0}
\end{bmatrix} - \sum_{(i,j)\in \mathcal{M}} y_{ij} M_{ij} - \sum_{(k,j)\in \mathcal{\bar{M}}} y_{kj} \bar{M}_{kj} \Big),
\label{eq:snl_sdr_gdam_init}
\end{equation}
where $\lambda $ is a positive number, which we choose to be $10.0$ in our implementation. After a few steps of GD, an feasible interior initialization for the dual problem \eqref{eq:snl_sdr_dual} can be readily obtained.

\subsection{Main solve}
After a strictly feasible initialization is found, the main solution phase starts that applies the accelerated Algorithm \ref{alg:2} for the dual problem \eqref{eq:snl_sdr_dual}. To further improve the practical performance, we developed adaptive stepsize and restart strategies, whose heuristics were tuned for the SNL problems. We implement a backtracking-like stepsize rule for the iterate to remain in the feasible set, i.e., the stepsize is reduced by a scaling factor whenever a constraint is violated. We choose the initial stepsize to be $\alpha_0 = 1.618e1$ and set the minimum stepsize to be $\alpha_{min} = 1e-8$. 
The restart strategy is motivated by the work \citep{o2015adaptive}. Specifically, we restart the optimization algorithm when the objective function value increases or the objective reduction is heavily slowed down\footnote{In our practical implementation, the objective function values are monitored in a fixed frequency.}. For efficiency, we do not reset the stepsize to $\alpha_0$, but to a multiple of the stepsize prior to the restart, i.e., $\alpha_{r+1} = \kappa \alpha_r$, where $\kappa$ is a constant and $r$ is the index of each restart. It should be noted that the implemented heuristic strategies are still basic, and more sophisticated methods can be used to further enhance the practical performance. Nonetheless, we found the present implementation is robust and scales well to large-scale problems. 
\subsection{Postsolve}
After solving the dual SDR problem \eqref{eq:snl_sdr_dual}, an additional step is needed to recover the primal solution $Z$ that contains the actual sensor locations. Denote the feasible sets of \eqref{eq:snl_sdr} and \eqref{eq:snl_sdr_dual} by $\mathcal{F}_p$ and $\mathcal{F}_d$, respectively. Denote $\mathcal{F} = \mathcal{F}_p \times \mathcal{F}_d$, and the interior of $\mathcal{F}$ by $ \overset{o}{\mathcal{F}}$. Assume that $\overset{o}{\mathcal{F}} \neq \emptyset$, the central path of a SDP can be expressed as
\begin{equation}
    \mathcal{C} = \{ (Z,y,S) \in \overset{o}{\mathcal{F}}: ZS = \eta I, 0 < \eta < \infty \},
\label{eq:sdp_cp}
\end{equation}
where $\eta$ is the barrier parameter (compare \eqref{eq:logarithmic_barrier}), $I$ is the identity matrix. The primal solution $Z$ can be recovered from the dual solution $S$ by
\begin{equation}
    Z = \eta S^{-1},
\end{equation}
where $\eta$ is the barrier parameter, and, within the computational framework of GDAM, it can be approximated by \eqref{eq:barrier_zeta},
\begin{equation*}
 \eta(x_\zeta^\sharp) =\zeta \frac{|\nabla f(x_\zeta^\sharp) |}{|\nabla \Phi(x_\zeta^\sharp) |},
\end{equation*}
where $x_\zeta^\sharp$ is a solution found by the parameter $\zeta$. To achieve higher accuracy, we run GDAM with $\zeta = 1.0$ by taking the optimizer of the main solve as the initialization. In our implementation, we set $\alpha_{min} = 1e-10$ for the postsolve. Note that smaller $\alpha_{min}$ (for both the main solve and the postsolve) leads to more accurate solutions. However, too small step sizes result in (almost) singular matrices for $S$, and numerical instabilities occur when inverting them.

\end{document}